\documentclass[a4paper,11pt,pointlessnumbers]{scrartcl}
\usepackage{amsmath}
\usepackage{amsfonts}
\usepackage{amssymb,mathrsfs}
\usepackage{t1enc}
\usepackage[latin1]{inputenc}
\usepackage[english]{babel}
\usepackage{natbib}
\usepackage{amssymb,stmaryrd}
\usepackage[final]{pdfpages}
\usepackage{caption}
\usepackage{url}

\usepackage[colorlinks=false, linktocpage=true,backref,pagebackref,hidelinks]{hyperref}

\newcommand{\Q}{\mathbb Q}
\newcommand{\C}{\mathbb C}
\newcommand{\R}{\mathbb R}

\newcommand{\M}{\mathbb M}

\newcommand{\E}{\mathbb E}

\newcommand{\As}{\mathscr{A}}

\newcommand{\Bs}{\mathscr{B}}
\newcommand{\Ds}{\mathscr{D}}
\newcommand{\Fs}{\mathscr{F}}
\newcommand{\Hs}{\mathscr{H}}
\newcommand{\Js}{\mathscr{J}}

\newcommand{\Ss}{\mathscr{S}}

\newcommand{\zbar}{\overline{z}}

\newcommand{\beq}{\begin{equation}}
\newcommand{\eeq}{\end{equation}}
\newcommand{\beqarr}{\begin{eqnarray}}
\newcommand{\eeqarr}{\end{eqnarray}}
\newcommand{\beqa}{\begin{eqnarray*}}
\newcommand{\eeqa}{\end{eqnarray*}}
\newtheorem{theorem}{Theorem}
\begin{document}
\thispagestyle{empty}

\title{\bf \large From Grassmann complements to Hodge-duality}
\author{\normalsize Erhard Scholz\footnote{University of Wuppertal, Faculty of  Math./Natural Sciences, and Interdisciplinary Centre for History and Philosophy of Science, \quad  scholz@math.uni-wuppertal.de}}
\date{\small 09/03/2020 } 
\maketitle
\begin{abstract}

\end{abstract}

\tableofcontents
\newpage

\vspace*{1em}
\section*{\small Introduction}
\addcontentsline{toc}{section}{\protect\numberline{}Introduction}

Hodge duality is a central concept of 20th century algebraic and analytic geometry  and plays a non-negligible role also in  recent mathematical physics. At first sight one might expect that its origins  lie in the 1930s when its name-giving protagonist, William V.D. Hodge, started his mathematical research. On the other hand, a close link between Hodge's theory and the Maxwell equation  has  been claimed not only from a systematic point of view but also historically  by Hodge' former student M. Atiyah in his talks on dualities in mathematics and physics.\footnote{``Maxwell's equations actually motivated Hodge for his work on harmonic forms in general. As indicated, Maxwell's equations are about forms of degree 2 in 4 dimensions and Hodge went to forms of any degree $q$ in any dimension $n$'' \citep[p. 77]{Atiyah:2008duality}. \label{fn Atiyah quote intro} }
 The question of how dense this connection was historically   leads back to the late 19th and early 20th century development of electrodynamics. Of course we are well advised not to take  systematic correspondences too easily as an indicator of a  historically effective relationships. Historical mis-readings of F. Klein's early attempt at reconstructing   Riemann's study of harmonic functions and meromorphic forms in complex dimension $n=1$ by linking it to the study of harmonic  flows  on surfaces \citep{Klein:RiemannscheFlaechen} may serve as a warning.

If one undertakes a journey to the late 19th and early 20th century with this question in mind,  an  unexpected author comes into sight: Hermann Grassmann.   Grassmann  was of some influence for the early understanding of the linear algebraic background for  duality concepts in electromagnetism and elsewhere. He and the  authors following him in this respect did not speak about ``duality'' but used the language of ``complements'' of alternating  products of  ``extensive'' (later vectorial) quantities. 
This may explain why  the thematic arc of this paper,   indicated already in the  title, is spanned  so wide. 

Our report starts with a glance at Grassmann's so-called complements of alternating products. Readers who are acquainted with  Hodge duality will immediately perceive it as a linear algebraic  template  (``precursor'') for the later Hodge $\ast$-operator (section 1). The first two subsections (1.1 and 1.2) discuss how Grassmann proceeded; subsection  1.3 puts the Grassmannian complements into the context of the vectorial operations developed simultaneously in the post-Hamilton tradition  independently  of  Grassmann's theory.  Readers who are not so well versed with Hodge duality and the Hodge  operator can  slowly approach  the topic by getting acquainted with this  more elementary twin. 

Section 2 jumps into the history of electromagnetism from an extremely selective viewpoint dictated by our topic. It does not try to trace all kinds of duality relations in electromagnetism, which early on in its history played an important heuristic role  and were often linked to philosophical, sometimes vague speculations on dualism between electricity and magnetism. Here we concentrate on that phase of electromagnetism in which duality in the sense of Grassmann complements became visible and was expressed in a clear mathematical form.  That means we jump into the history of electrodynamics when it became relativistic. Readers who want to get a broader perspective of electromagnetism and even only the Maxwell equation  have to consult the respective literature in the history of physics.\footnote{See in particular  \citep{Darrigol:Electrodynamics,Hunt:1991,Steinle:2016,Steinle:Feldtheorie}.}
 Like so many fundamental concepts in special relativity the introduction of Grassmann complements into the formulation of the Maxwell equations stems from Hermann Minkowski (sec. 2.1). The linear algebraic duality relation contained in Minkowski's relativistic electrodynamics  was taken over into the more involved ``curved'' background of general relativity, i.e. was adapted to the infinitesimal structures (the tangent and cotangent spaces/bundles and their sections, in more recent terminology) of Riemannian geometry. Several authors played here a role, among them notably  W. Pauli, F. Kottler, and T. de Donder (sec. 2.2). Einstein did not belong to this group; he preferred to formulate Maxwellian electrodynamics without an (explicit) reference to this type of duality relation.

  Elie Cartan, on the other hand, did so (sec. 2.3). It may come as a surprise that  Cartan even used a generalized version of Grassmann complements not only for expressing  the Maxwell equation, but also for his reformulation and generalization of Einstein gravity. In his attempt at  bridge-building between Einstein's theory and the generalized theory of elasticity of E. and F. Cosserat Cartan used a transmutation of translational curvature into one represented as a rotational momentum (torsion) and vice versa.\footnote{See the separate and more detailed discussion in \citep{Scholz:2019Cartan}.} 
In three-dimensional space, which served Cartan as a guide to his intuition, a  Grassmann type duality was just the right tool; Cartan  then used it also for  formulating   his generalized theory of gravity in dimension $n=4$.

The third section  discusses the origins of Hodge duality. In order to do justice to Hodge, a short introduction to the background of Riemann's theory of Abelian integrals is necessary (sec. 3.1). Hodge's research started with the goal to do the ``same'' as Riemann (or something similar)  for $n$ dimensions, which required as a first step to define what it means for a differential form on a manifold $M$ with metric (a Riemannian manifold) to be ``harmonic''. A short outline of the early stage of this research program is given in sec. 3.2.  Hodge's research soon let him state an intriguing relationship, although not immediately proven in a valid form, between the numbers of linearly independent harmonic $p$-forms and the topology of $M$, viz. its Betti numbers. This led him to new insights into the theory of algebraic surfaces  (sec. 3.3). 

In the mid 1930s  and 1940s Hodge brought his program into  a first stage of maturity, documented in his book on {\em The Theory and Application of Harmonic Integrals} \citep{Hodge:1941}. Here  all the basic topics  we are interested in can be  found  in sufficient clarity: Hodge $\ast$-operation, the definition of harmonic forms, and  a discussion of the relation to the Maxwell equation (sec. 3.4). It may be disillusioning, however, to see that  his discussion of the Maxwell equation remained quite  elementary and reduced to simple special cases. That  seems to be characteristic for Hodge's research profile; he was rather a (pure) mathematician than truly interested in the link to mathematical physics.

 Subsection 3.5 gives a condensed picture of how Hodge's methods and results were assimilated by the wider community of complex algebraic geometers as far as transcendent (i.e. analytic) methods were concerned.  This section should not be misunderstood  as a full history of this assimilation; it can only give short glimpses into the phase in which Hodge's theory  was generalized and turned into what since then became to be known as {\em Hodge theory}, with the Hodge-Laplacian  and  the study of the Hodge structure of holomorphic and anti-holomoprphic forms on a complex K\"ahlerian manifold.

Developments after mid-century  are only touched at in the outlook of section 4. It would be inadequate, however,to stop short in our report without at least mentioning the  shift Hodge theory underwent with the introduction of sheaf cohomological methods by P. Dolbeault,  J.-P. Serre and others. In the end, Serre's famous duality theorem  generalized  Hodge duality and had lasting consequences for  mathematics  in the second half of of the twentieth century (sec. 4.1). The final subsection 4.2    turns  back to the question of how Hodge's duality concept was taken up in  physics. The strongest role for it  arose  within  Yang-Mills theory. 

In  the 1970s and 1980s, the decades of  the rising standard model of elementary particles, the manifest success for Young-Mills type gauge fields in physics had  strong repercussions also in mathematics. It inspired mathematicians to study the vacuum solutions of Yang-Mills equations. This lead to a kind of generalized theory of harmonic forms and became a  fruitful tool in differential topology, in particular for manifolds of dimension $n=4$.  The remarks of our outlook on this topic  remain of course incomplete. They cannot be more than first hints to a field of recent mathematics the history  of which has still to be written.

\section{\small Grassmann's complement of ``extensive quantities'' \label{section Grassmann}}

\subsection{\small Grassmann's {\em Ausdehnungslehre} (1862) \label{subsection Ausdehnungslehre}}
{\em Hermann G\"unther Grassmann} (1809--1877) is today often recalled as the founder of $n$-dimensional vector spaces and the algebra of their exterior products, clad in a peculiar  geometrico-algebraic theory which he called  the {\em doctrine of extension} ({``Ausdehnungslehre''}).\footnote{See among others \citep{Petsche:GrassmannEnglish,Schubring:Grassmann}}
 He got  acquainted with vector-like  ideas through his father's, {\em Justus G\"unther Grassmann}'s,  publications on  crystallography  and gave a condensed review of them \citep{Grassmann:1839} before he embarked on the project of his general theory of extension.\footnote{\citep{Scholz:Justus-Hermann-Grassmann} }
  He published the  ideas on his   {\em Ausdehnungslehre}  in two book-length general publications \citep{Grassmann:1844,Grassmann:1862}, the second one now also available in English \citep{Grassmann:1862English}, and several articles  on more   specialized topics.\footnote{For the articles see \citep[vols. 2.1, 2.2, 3.1]{Grassmann:Werke}. The first {\em Ausdehnungslehre} (1844) is now also available in French \citep{Grassmann:Francaise}.}  
 
For present day  readers  the exposition in the second book is  easier to follow than the first one. Here Grassmann  introduced an {\em extensive quantity} (``Ausdehnungsgr\"osse'') $a$ of an {extensive domain of level $n$} (``Ausdehnungsgebiet $n$-ter Stufe'') as a linear combination of a finite number of independent generators $e_1, \ldots , e_n$ (Grassmann's symbols), called by him the  {\em original units} (``urspr\"ungliche Einheiten'') 
 \beq a = \sum_{i=1}^n \lambda_i e_i \, ,
 \eeq
 where  $\lambda_i$ denote real, rational or later even complex coefficients. On this basis he  explained the  operations  (addition, scalar multiplication) and their fundamental laws \citep[chap. 1, part I]{Grassmann:1862}, which were later condensed  into the axioms of a vector space  \citep{Peano:Calcolo_geometrico}.\footnote{\citep{Dorier:Vector}}
 In this way a Grassmannian extensive domain of level $n$ may  essentially be understood  as an {\em $n$-dimensional vector space} $V$ over the rationals, the real or the complex numbers.  
 
 He  discussed bi- and multilinear products of two or more extensive quantities  $a, b, \ldots$ in general (chap. 2) -- in later terminology {\em tensor products} -- and showed  particular interest for alternating products of  quantities written by him in square brackets (without commas)  $[ab]$ or, in the case of several, say $k$,  independent quantities $ a_1, \ldots, a_k $, $(1\leq k\leq n)$   as
 \[
 A= [a_1 \ldots a_k] \, .
 \]
 Grassmann  called this a {\em combinatorial product} or {\em exterior product} (``\"ausseres Produkt'')  (chap. 3, part I). 
Such products  formed   {\em quantities of the $k$-th level}  (``Gr\"ossen  $k$-ter Stufe'')  and were  interpreted by him geometrically -- we would say as  an equivalence class of parallelotopes of dimension $k$. These quantities were thus of a  new type; by linear combination they established a new {\em extensive domain of $k$-th level} in a {\em principal domain of level $n$}  (``Hauptgebiet $n$-ter Stufe''). They could be generated by ``simple'' quantities of $k$-th level of form $[e_{i_1} \ldots e_{i_k}]$ with original units $e_{i_1},  \ldots, e_{i_k}$. 
  In modernized symbolism such an extensive domain of $k$-th level corresponded to the $k$-th exterior product $\Lambda ^k V $ of $V$
with elements $A$ linearly generated by quantitites of the type $[a_1 \ldots a_k]$.  
 
Of course Grassmann realized that for $k=n$ the domain of the $n$-th level is generated by the combinatorial product of all original units. He therefore  identified $\Lambda ^n V$ with he ``numerical quantities''  $\mathfrak{N}$ by the following rule \citep[\S 89]{Grassmann:1844}\footnote{``{\em Erkl\"arung}: Wenn in einem Hauptgebiete $n$-ter Stufe das kombinatorische Produkt der urspr\"uglichen Einheiten $e_1, e_2, \ldots, e_n$ gleich 1 gesetzt ist und $E$ eine Einheit beliebiger Stufe, (\ldots), so nenne ich ``{\em Erg\"anzung} von $E$ \ldots ''  \citep[\S 89]{Grassmann:1844} (cf. next footnote).  
Today one would substitute the number field $\Q, \R$ or $\C$ for   $\mathfrak{N} $.}
\begin{quote}
{\bf Unit convention}:  The combinatorial product of all original units is considered to be  equal to the ``numerical unit'':
\beq [e_1 e_2 \ldots e_n ] =1 \; 
\eeq 
\end{quote}
On such a basis he studied the properties of the alternating algebra over $V$  in quite some detail, e.g., the rules of the operation which we would  write as $\Lambda ^k V \times \Lambda ^l V \longrightarrow \Lambda ^{k+l} V$. He did so  even for the case $k+l > n$ by reducing the result ``modulo $\Lambda ^n V$'' and the unit convention. Today such a structure is  called a {\em Grassmann algebra}. 

Only in the next step of developing his theory (chap. 4, part I) Grassmann introduced an {\em interior product} (``inneres Produkt'') in an extensive domain (see below). He did so in such a way that the original units $e_1, \ldots, e_n$ turned into  an {\em orthonormal basis} of the  extensive domain, the vector space $V$, to  use  later terminology. The introduction was built upon the properties of the  complement (``Erg\"anzung)'' of an extensive quantity. This is not all we find in \citep{Grassmann:1862}. In chapter 5 (part I)  of the book the author sketched how to  study the  classical geometrical, i.e. Euclidean, space as an extensive domain of level 3 with exterior/combinatorial and interior products. In  part II he studied what we would call  analysis in $n$-dimensional real vector spaces with differentiation, infinite series and integration.   

We need not look at  these interesting parts of his book here; for our purpose it suffices to concentrate on the introduction of Grassmann's (dual) complement in an extensive domain.

\subsection{\small Dual complement ({\em Erg\"anzung})\label{subsection dual complements}}
In \S 4 of his chapter 3  on combinatorial products Grassmann introduced a new  operation by forming the {\em complement} of a given quantity (``Erg\"anzung der Gr\"ossen'', article 89). For $E$ a {\em unit of the $k$-th level}, in slightly modernized notation (using double subscripts) $E = [e_{i_1} e_{i_2} \ldots e_{i_k}]$, he introduced the  
\begin{quote}
{\bf Definition}: The {\em complement} $E'$ of a unit of the $k$-th level $E = [e_{i_1} e_{i_2} \ldots e_{i_k}]$ is the combinatorial product of all units not appearing in $E$, say $E'= \pm [e_{j_1}e_{j_2} \ldots e_{j_{n-k}}]$, where the sign is chosen in such a way that\footnote{``\ldots nenne ich `Erg\"anzung von $E$' diejenige Gr\"osse  $E'$, welche dem kombinatorischen Produkte $E'$ aller in $E$ nicht vorkommenden Einheiten gleich oder entgegengesetzt ist, je nachdem $[EE']$ der absoluten Einheit gleich oder entegegengesetzt ist.'' \citep[p. 57, art. 89]{Grassmann:1862} }  
\[  [E E']= [e_1 e_2 \ldots e_n] \quad \mbox{is the ``absolute unit''} \, .
\]
Grassmann used the standard notation  $|E $ for the complement of $E$. 

For any quantity $A$ of $k$-th level, say $A =  \alpha_1 E_1 + \ldots \alpha_l E_l$  (with $E_1$, \ldots, $E_l$ units of the $k$-the level) the complement is defined by  linear continuation:
\[ |A =  \alpha_1 | E_1 + \ldots \alpha_l | E_l
\]  
\end{quote}
Using the later notation $\wedge$ for Grassmann's combinatorial/exterior  products,  the complement of $E = e_{i_1} \wedge e_{i_2} \wedge   \ldots \wedge e_{i_k}$ (where the $e_i$-s are pairwise distinct) is given by
\beq |E= \mathsf{sg}(i_1, \ldots i_k, j_1, \ldots j_{n-k})\, e_{j_1}\wedge \ldots \wedge e_{j_{n-k}} \, 
\eeq 
where $\{{i_1}, \ldots , {i_k},  {j_1},  \ldots {j_{n-k}} \} =\{ 1, \ldots, n\}$ and $\mathsf{sg}(\ldots)$ denotes the sign of the permutation $(\ldots)$.

Grassmann studied   the operation
\[  |: \; \wedge^k V \longrightarrow \wedge^{n-k}V
\] 
in quite some detail. Among others  he showed:
\beq || A  =  (-1)^{k(n-k)} A \qquad \mbox{\citep[art. 92]{Grassmann:1862}}
\eeq
In particular for $n$ uneven $||A=A$, while for $n$ even $||A = (-1)^k A$. Up to sign, the complement of the Grassmann complement is the original quantity. After the turn to the 20th century this would be considered  as a kind of duality and $|A$ was accordingly called the dual complement (``duale Erg\"anzung'' of $A$), for example in \citep[\S 12]{Pauli:1921}; see section \ref{subsection GTR}.

For the Grassmann complement of a product, $A\in\wedge^k V , \, B \in  \wedge^l V$  (and  $\mathit{dim}V = n$) our author showed
\beq |\, [AB] = [\, |A\, |B]  \label{dualization of products}
\eeq 
if $k+l < n $ \citep[art. 94]{Grassmann:1862}. Using his definition of the combinatorial product for $k+l>n$ by factorizing out products $[e_1 \ldots e_n]\equiv 1$ Grassmann  extended the products and the complements to the general case. In this constellation he spoke of a ``regressive product'', in contrast to a ``progressive product'' in the case $k+l\leq n$,  and found that   (\ref{dualization of products}) can be extended to the general case \citep[art. 97]{Grassmann:1862} and to an arbitrary number of factors (art. 98). In particular for extensive quantites $a, b \ldots f$ of the first level (for us, elements of $V$)
\beq | \, [ab \ldots f] = [ \, |a |b \ldots |f] \, .
\eeq

Using the dual complement and the exterior product  Grassmann finally defined the {\em interior product} of two extensive quantities $A, B$ (of any level) as 
\beq [A\, |B] \,  \qquad \mbox{\citep[art. 138]{Grassmann:1862}.}\eeq
The level of the interior product is $n+ k -l$ for $l>k$ and $k-l$ otherwise. The {\em interior product of quantities of the same level is thus a numerical quantity}
 (``Zahlgr\"osse'', art. 141). For  $k=l=1$ and $n=3$ 
it coincided with Hamilton's {\em scalar product} of vectors. For general dimension $n$  Grassmann had  his version of what later would be called a positive definite bilinear form, or Euclidean scalar product.

\subsection{\small Vector operations in $\R^3$ in a Grassmann's perspective \label{subsection vector operations}}
Grassmann's theory of extensive quantities was not immediately absorbed by mathematicians and/or physicists. As a symbolical tool for mathematical physics the vectorial calculus in three-dimensional space arising from Hamilton's quaternions was a strong competitor during the last third of  the 19th century.\footnote{\citep{Crowe:Vectors,Reich:Grassmann}}
 But the reception of Grassmann's theory was not exceptionally slow if one compares it with other general mathematical theories  \citep{Rowe:Grassmann,Tobies:Grassmann}.\footnote{See, in particular  \citep[p. 132]{Rowe:Grassmann}.}

 Early on, i.e. already in the years 1844ff., also Hamilton  insisted on the usefulness of the non-commutative product of the vector part of quaternions, inherited from the product structure of the quaternions. What he called  the {\em vector product} stood in close relation to Grassmann's complement of the exterior product of two extensive quantities in 3-space. His (commutative) {\em scalar product} was equivalent to a special case of Grassmann's interior product.   At the turn to the 1880s Josia W. Gibbs was one of the influential authors who separated the 3-dimensional vector calculus from the quaternionic calculus. The other one was Oliver Heaviside, see \citep[chap. 5]{Crowe:Vectors}.  
Gibbs knew Grassmann's approach well and  was aware of the relation  between the vector product in the sense of Hamilton and the (dual complement of the) exterior product for three-dimensional ``extensive quantities''.\footnote{In the draft of a letter to Schlegel in 1888 Gibbs wrote: ``My acquaintance with Grassmann's work was also due to the subject of E. [electricity] \& in particular to the note wh (sic!) he published in Crelle's Jour, in 1877 calling attention to the fact that the  law of the mutual action of two elements of currents  wh Clausius had just published had been given  1845 by himself. \ldots  that law is so very simply expressed by
means of the external product'' cited from \citep[p. 153]{Crowe:Vectors}. Gibbs'  referred her to the original publication  \citep{Grassmann:1845}. } 

Here the special type of Grassmann duality in 3-dimensional   space, in short $V\cong \R^3$  endowed with the Euclidean interior product and ``units'' (orthonormal basis vectors) $e_1, e_2, e_3$  came into play. 
  In terms of  Grassmann's exterior/combinatorial product  and his dual complement the vector product $w$ of $u,v \in V$  is nothing but 
  \[ w = u \times v  = |\, [u\,v]
  \] 
  Grassmann knew of course  the relation -- valid in all dimensions --  which we would write as mappings
  \beq |: \;  \wedge^{n-1}V \stackrel{\cong}{\longrightarrow} V\, \qquad |:\;  [e_1 \ldots e_{n-1}] \mapsto e_n \quad \mbox{etc.} 
  \eeq
  (with etc. indicating that cyclical permutations of indices are permitted)
   and could easily use it for introducing an alternating product in 3-dimensions. For an orthonormal basis this implies in particular the identity $e_1 \times e_2= | \, [e_1 e_2] = e_3$ etc. 
   
   In the last third of the 19th century a varied scene of vector algebra and analysis developed among mathematicians and
    physicists.\footnote{\citep[chaps. 5--7]{Crowe:Vectors}, \citep{Reich:Grassmann}.}
 In this  heterogeneous milieu  Grassmann's exterior product of vectors (``extensive quantities of 2-nd level'') in 3-space were often assimilated to vectors, often suppressing the duality operation $|$ mediating between the two types of quantities. But also in the physics tradition a distinction between the two was realized to be important. The reason was that  point inversions, which became important in crystallography, have different effects on elements in $V\cong \R^3$ and in $\wedge^2 V$. In the first case an extensive quantity changes sign, in the second case it does not (corresponding to the  eigenvalues  of an inversion $-1$ respectively $+1$ in later terminology).

 In his {\em Kompendium der theoretischen Physik} (1895/96) {\em  Woldemar Voigt} introduced the terminology of a {\em polar vector} for an  element in $V\cong \R^3$ and {\em axial vector} for one  arising as the complement of an external product, i.e. an element of  $| \wedge^2 V $ (where, of course, the notation mixing Grassmann's symbol $|$ with a modern one for the alternating product is mine).\footnote{See \citep[p. 10, fn. 17]{Abraham:1901} and  \citep[p. 203]{Reich:Grassmann}.} 
  From a physical point of view the most important elementary case of an axial vector was a ``couple of vectors'' $(u,v)$ acting on a rigid body. The  combined action of a couple (a rotational momentum in later terminology) could best be symbolized as $u\times v$ or $ |\, [u v]$, respectively.
This terminology entered the survey article on ``Basic concepts of geometry''  in {\em Enzyklop\"adie der Mathematischen Wissenschaften}, written by Max Abraham \citep[p. 6ff.]{Abraham:1901}. It may thus  be considered as part of the general knowledge at the turn to the 20th century.


\section{\small Dual complements in physics of the early 20th century \label{section dual complements in physics}}

\subsection{\small Background: Maxwellian electrodynamics in the 19th century}
In mechanics and the geometry of the $\R^3$ Grassmann type duality could be handled   more easily in the framework of the 3-dimensional vector calculus arising from the  Hamilton tradition. This could be done without even noticing the underlying duality relationship which was well hidden in the rules of  calculation. Duality ideas appeared  in physics more directly in the context of electromagnetism,  initially in the natural philosophical guise of a  conjectured ``duality'' between electricity and magnetism, later extended to a kind of double pairing between the field strengths and field excitations of both, electricity and magnetism. Here we jump over the long history of two thirds of the 19th century\footnote{See \citep{Steinle:2013,Steinle:2016,Darrigol:Electrodynamics}.}
 and pass directly to 
the field theoretic formulation of electric and magnetic phenomena due mainly to {\em Michael Faraday } and {\em James Clerk Maxwell}. This  was a major achievement of physics in the 19th century. 

 As we have  already seen,  Gibbs' interest in Grassmann was partially due to the latter's first steps in the direction of a vector analytic representation of magnetic and electrical phenomena. In the study of electromagnetism certain aspects of  ``duality'' in the wider sense came up. At the time of Oerstedt and Faraday  electric and magnetic phenomena were considered as a striking, and heuristically successful example for the search of ``duality'' in sense of natural philosophy. Here it  found an expression in the pairing of the electric  field strength $\mathbf{E}$ (``electromotive force'') with the magnetic field $\mathbf{B}$ (``magnetic induction'' mentally depicted by Maxwell as a system of directed vortices formed by the magnetic fluid)  and their  excitations $\mathbf{D}$ (``electric displacement'')  and $\mathbf{H}$ (``magnetic intensity''). All of them were directed quantities in space, which slowly came to be known as {\em fields}.\footnote{\citep{McMullin:fields,Steinle:Feldtheorie}}

Close to the end of the 19th century   the (non-relativistic) Maxwell equations were written in different notations (in components, quaternionic, or vector analysis) and with changing conventions for the proportionality constants. In slightly modified notation (for the field strengths and excitations) and using  3-dimensional vector calculus they can be resumed  in a form  similar to the one used in \citep{Hertz:1890}, and \citep{Foeppl:Maxwell} (which Einstein knew well):\footnote{Similarly in papers by O. Heaviside in the second half of the 1880s. For details see  \citep[p. 125f.]{Hunt:1991} and \citep{Darrigol:Electrodynamics}.}
\beqarr \mathsf{curl}\, \mathbf{E} + \frac{1}{c}\partial_t \mathbf{B} &=& 0 \, , \qquad \mathsf{div}\, \mathbf{B} = 0 \label{eq Maxwell I class}\\
\mathsf{curl}\, {\mathbf{H}} -\frac{1}{c} \partial_t \mathbf{D}s &=& \mathbf{J} \; , \qquad \mathsf{div}\, \mathbf{D} = \rho \, , \label{eq Maxwell II class}
\eeqarr
where $\rho$ and $\mathbf{J}=(J_1,J_2,J_3)$ denote the charge and current densities respectively and $c$ the velocity of light. The two equations of (\ref{eq Maxwell I class}) are a refined mathematical expression for Faraday's induction law and the non-existence of magnetic charges  \citep[p. 449ff.]{Steinle:Feldtheorie}, while  (\ref{eq Maxwell II class})  can be read as a dynamicized version of Oerstedt's and Amp\`ere's characterization of the magnetic effects of an electric current taking Maxwell's displacement current into accout (ibid., p. 453ff.), and the existence of charges as the source of an electric displacement field.  
The  quantities $(\mathbf{E},\mathbf{H})$, usually considered as the electric respectively magnetic field strength,  and their respective  excitations $(\mathbf{D},\mathbf{B})$ are assumed to be proportional,\footnote{In the literature a variety of  terminology is to be found; often $\mathbf{D}$ is called the electric displacement, $\mathbf{B}$ the magnetic induction. Here both are being considered to be excitations of the corresponding field strengths.} 
  with  (material dependent) coefficients the  dielectric constant $\epsilon$ and the magnetic permeability $\mu$,  such that
\beq \mathbf{D} =  \epsilon\, \mathbf{E}   \, , \qquad    \mathbf{B} = \mu\,  \mathbf{H} \, . \label{eq constants}
\eeq 
Moreover $\epsilon_0 \mu_0=c^{-2}$ for the vacuum values  $\epsilon_0, \mu_0$ of the coefficients. 

This  way of writing the four Maxwell equations  suggests a close analogy  between the pairs  $(\mathbf{E},        \mathbf{B})$ and  $(\mathbf{D},    \mathbf{H})$, expressed in the two sets  (\ref{eq Maxwell I class})  and (\ref{eq Maxwell II class}) of the Maxwell equations. The analogy seemed  far away from Grassmann duality, but that changed with the relativistic reformulation and generalization of Maxwellian electrodynamics. The step from the non-relativistic to the relativistic theory of the Maxwell-Lorentz theory is a highly involved story which deserves more interest than can be invested here; in our discussion we focus on an extremely selective aspect, the introduction of duality concepts in the sense of Grassmann. The following investigation  must therefore neither be understood as a history of electrodynamics, nor does it claim to do justice to the complicated  history of the linear algebra in the early part of the 20th century, which would make it necessary to relate it to the influence of the Berlin tradition of matrices put forward by L. Kronecker, K. Weierstrass and G. Frobenius.\footnote{See \citep{Hawkins:Weierstrass,Hawkins:2005,Brechenmacher:2016}}
 The selection of texts discussed in the next sections is governed by the criterion to trace the reappearance of Grassmann type duality ideas, sometimes in explicit reference to the latter and sometimes even without knowing about the Grassmannian background. 

\subsection{\small Einstein's relativistic electrodynamics,  Minkowski's ``six-vectors'' and their duals \label{subsection Electrodynamics}}

We enter this story with a short glance at {\em Albert Einstein}'s 
conceptual analysis of the Lorentzian relativity principle in his famous paper on  ``The electrodynamics of moving bodies'' \citep{Einstein:1905SRT}. There he considered, among others, the transformation of an  electric field $\mathbf{E}=(E_1,E_2,E_3)$ and a magnetic field given by $\mathbf{B}=(B_1,B_2,B_3)$ from one system of reference $\mathfrak{S}$ to another  one  $\mathfrak{S}'$ which is in linear uniform motion with respect to $\mathfrak{S}$. An electric charge $q$ resting in  $\mathfrak{S}'$ will be in motion, if considered from the point of view of  $\mathfrak{S}$. This leads to  magnetic forces induced by the motion of $q$ with regard to $\mathbf{E}$ in addition to those of $\mathbf{B}$. Einstein considered the transformation laws for the combined field quantities $(\mathbf{E},\mathbf{B})$ and concluded that 
\begin{quote}
\ldots the magnetic and electric forces have no existence independent of the state of motion of the coordinate system.\footnote{``\ldots da\ss{} die elektrischen und magnetischen Kr\"afte keine von dem Bewegungszustande des Koordinatensystems unabh\"angige Existenz besitzen''  \citep[p. 910]{Einstein:1905SRT}.}
\end{quote}

Einstein thus indicated an underlying unity of the fields $(\mathbf{E},\mathbf{B})$, and similarly  ($\mathbf{D},\mathbf{H}$),  but he was not able to cast the intended unification into  a proper mathematical form beyond the transformation formulas. This  was achieved  a few years later by {\em Hermann Minkowski}. 

Minkowski proposed a  unified mathematical representation of the (special) relativistic electromagnetic field by two antisymmetric matrices $F=(F_{ij})$ and $f=(f_{ij})$.  
 The first one $F=(F_{ij})$,  later on called the {\em Faraday tensor},  
 contained  the 6 components of  $(\mathbf{E},  \mathbf{B})$ in a proper arrangement,  and  $f$    combined the components  of $(\mathbf{D}, \mathbf{H})$, where $F_{ij}=-F_{ji}$ and  $f_{ij}=-f_{ji}$ \citep[pp. 356, 168]{Minkowski:1908Grundgleichungen}.\footnote{$F_{12}=B_3, F_{13}=-B_2, F_{23}=B_1$ and  $ F_{j4}=-\sqrt{-1}\, E_j,\; j=1, \ldots 3 $;  similarly for $f$ with $B_j, E_j$ replaced with $H_j$  and $D_j$ respectively. The factor $-\sqrt{-1}$ is due to Minkowski's notation of the Lorentz metric with imaginar time components. With the exception of $E$ Minkowski used different symbols,  $\mathfrak{m}=\Hs, \mathfrak{M}=\Bs$  and ${ \bf e}=\Ds$. \label{fn F Minkowski}}
  Moreover, Minkowski used units such that $c=1$  (at least for the numerical value). Then $\epsilon_0= \mu_0^{-1}$. and the  two matrices were  proportional,
  \beq f= \mu_0^{-1}F\, . \label{vacuum proportionality}
  \eeq
 
He explained how these matrices operate as antisymmetric bilinear forms on vectors in spacetime (which later became to be known as {\em Minkowski space})  $x=(x_1, \ldots, x_4)$ and $u=(u_1, \ldots, u_4)$, and how they have to be modified under Lorentz transformations. In other words he used the matrices to represent antisymmetric tensors  or alternating 2-forms  over  spacetime.\footnote{In short notation $f(x,u)= ^tx \cdot  f \cdot u $ and $F(x,u)= ^tx \cdot  F \cdot u $ \citep[eq. (23) (24), p. 364f.]{Minkowski:1908Grundgleichungen}.}
Minkowski neither used the terminology of tensors, nor did he mention Grassmann's alternating  products. He rather talked about his matrices as  {\em spacetime vectors of the second kind}. This language  was much closer to the post-Hamilton vectorial calculus in 3-dimensional space  than to  Grassmann,  Here ``axial vectors'' like   angular velocity and rotational momenta were also called vectors of the second kind. But in 4-dimensional space the identification with the original vector space was no longer feasible and Minkowski had to introduce more general  linear algebraic methods which may  have been inspired by  Berlin tradition of matrices and representations.\footnote{I owe this observation to a hint by F. Brechenmacher.}
 Anyway, for the next few years Minkowski's antisymmetric tensors  became to be known as ``six-vectors" (Sommerfeld's terminology),\footnote{\citep{Walter:Breaking_in}}  before the tensor terminology took roots with the general theory of relativity.

Minkowski introduced a symbolism of his own for the vector analysis on Minkowski space, extending , e.g., the $\nabla$ operator of 3-dimensional vector analysis to an operator called $\mathsf{lor}=(\partial_1, \ldots , \partial_4)$ by him, which transformed under Lorentz transformations like a vector in Minkowski space. 
For the second set of the  Maxwell equations (\ref{eq Maxwell II class}) he  got, in slightly adapted notation, 
\beq \mathsf{div}\, f =   s\,, \label{eq Maxwell II Minkowski}
\eeq 
in the sense of $\sum_j \partial_j f_{ij}=  s_i$, where $s$ denotes the 4-current with $s_k=J_k$ for $j=1,2,3$ and, essentially,\footnote{I delete  the imaginary factor $\sqrt{-1}$ used by Minkowski for timelike components in order to formally assimilate the Minkowski metric to a Euclidean signature.} $s_4=\rho$. Moreover we have 
to keep in mind that he wrote the equation as $\mathsf{lor}\, f =  \sum_i \partial_i f_{ij}= - s_i$ \citep[p. 384]{Minkowski:1908Grundgleichungen}.

In order to bring the first set of the Maxwell equations (\ref{eq Maxwell I class}) into a comparable form, he introduced the {\em dual matrix} (``duale Matrix'')  $f^{\ast}$ (Minkowski's terminology and symbolism) of a skew symmetric matrix $f$ as
\beq f^{\ast}= (f^{\ast}_{ij})\; ,  \qquad \mbox{where} \quad f^{\ast}_{ij}= \mathsf{sg}(\mbox{\footnotesize ijkl}) f_{kl}\; 
\eeq 
with $\mathsf{sg}(\mbox{\footnotesize ijkl})$  the sign of the permutation. Then the first set of the Maxwell equation became in adapted notation
\beq \mathsf{div}\, F^{\ast} = 0 \;,  \label{eq Maxwell I Minkoswki}
\eeq 
which in Minkowki's own notation of 4-dimensional vector analysis read as $\mathsf{lor}\, F^{\ast} =0$. In the result Minkowski had shown how his new 4-dimensional  symbolism could be used to give a unified expression of the Maxwell equations. In order to achieve this he introduced an apparently new type of duality $A \mapsto A^{\ast}$ for skew symmetric matrices $A$.

 But how new was this type of duality? Seen from the perspective of Grassmann's  theory this duality was already inherent in the  complement of combinatorial products. One only needed to read Minkowski's 
matrices  as a pragmatic notation for  combinatorial  products of 1-forms over Minkowski space considered as an ``extensive domain'' generated by $e_1, \ldots, e_4$ (with respect to an orthonormal reference system $\mathfrak{S}$). Denoting the 1-forms as $e_i^{\ast}$ ($1\leq i \leq 4$)) with $e_i^{\ast}(e_j)=\delta^j_i$, the matrices  stand for
\beq F = \sum_{1\leq i<j \leq 4}F_{ij}[e_i^{\ast} e_j^{\ast}] \, , \qquad  f = \sum_{1\leq i<j \leq 4}f_{ij}[e_i^{\ast} e_j^{\ast}] \, ,
\eeq 
 in the more recent Cartan symbolism  of  differential forms   $F=\sum F_{ij}\,dx^i\wedge dx^j $ and $f=\sum f_{ij}\,dx^i\wedge dx^j $. Minkowski  did not mention Grassmann while he mentioned Hamilton's quaternions, although he did not find this calculus useful for his purpose \citep[p. 375, fn.]{Minkowski:1908Grundgleichungen}. This seems to indicate that Minkowski ws not aware of the possibility to express his dualization of matrices in Grassmannian terms. 

This was different for {\em Arnold Sommerfeld}, an important early promoter of Minkowski's approach to special relativity \citep{Walter:2010}. He  explained carefully and with an explicit reference to Grassmann how one could perceive  Minkowski's vectors of the second kind geometrically as a linear combination of ``surface elements of content 1''  in 4-dimensional spacetime \citep[p. 750]{Sommerfeld:1910SRT}. As they  have the dimension 6 he introduced the terminology of {six-vectors}, i.e. elements of $\Lambda^2 \M$ in our notation (where $\M$ denotes the Minkowski space with coordinates $(x_1, \ldots, x_4)$).\footnote{According to Sommerfeld, Emil Wiechert had already characterized the magnetic excitation $H$ (here called ``magnetische Feldst\"arke'')  as a Grassmannian quantity of the second kind in classical space \citep[p. 750]{Sommerfeld:1910SRT}.}  Projecting such a ``surface element of content 1'' $\varphi$ into the $(x_i,x_j)$-plane  gave him the coordinates $\varphi_{ij}$. That allowed him to define the {\em complement} (``Erg\"anzung'') $\varphi^{\ast}$ in agreement with Grassmann and with Minkowski by  \citep[p. 756]{Sommerfeld:1910SRT}
\[ \varphi^{\ast}_{ij}= \varphi_{kl}\; \qquad \mbox{with} \qquad \mathsf{sg}(ijkl)=1 \, .
\]
In this way Sommerfeld geometrized Minkowski's algebraic (matrix) approach to the Maxwell field, brought it into an explicit relation with Grassmann's concepts,  and prepared the way for a generalization to a varying (``curved'') metric on spacetime.

\subsection{\small Dual complements in the general theory of relativity \label{subsection GTR}}

Minkowski/Sommerfeld's   dualization  did not find   unanimous approval among physicists. Einstein even reproached such a  representation of the Maxwell field quantities as  ``involved and confusing''. In the introduction  to his paper on electrodynamics in the general theory of relativity  he  described  Minkowski's characterization of the electrodynamic field by a ``six-vector''   and  a second  dual  one `` whose components have (\ldots) the same values as the first one, but are distinct in the way the components are associated with the four coordinate axes'' \citep{Einstein:1916Maxwell} in a pragmatic way.   He correctly continued that then the
\begin{quote}
\ldots two systems of Maxwellian equations are obtained by setting the divergence of the first one equal to zero, and the divergence of the other one equal to the four-vector of the electric current.\footnote{``\ldots Man erh\"alt die beiden Maxwellschen Gleichungssysteme, indem man die Divergenz des einen Sechservektors gleich Null, den andern gleich dem Vierervektor des elektrischen Stromes setzt'' \citep[p. 184]{Einstein:1916Maxwell}.}
\end{quote}
But he was dissatisfied with such a characterization. He went on criticizing:
\begin{quote}
The introduction of the dual six-vector makes its covariance-theoretical representation relatively involved and confusing. Especially the derivations of the conservation theorems of momentum and energy are complicated, particularly in the case of the general theory of relativity, because it also considers the influence of the gravitational field upon the electromagnetic field (ibid.).\footnote{``Die Einf\"uhrung des dualen Sechservektors bringt es mit sich, da\ss{} diese kovariante Darstellung verh\"altnism\"a\ss{}ig  un\"ubersichtlich ist. Insbesondere gestaltet sich die Ableitung des Erhaltungssatzes des Impulses und der Energie kompliziert, besonders im Falle der allgemeinen Relativit\"atstheorie, welche den Einflu\ss{} des Gravitationsfeldes und des elektromagnetischen Feldes mitber\"ucksichtigt.'' (ibid.)}
\end{quote} 

In the light of later development, in particular Cartan's calculus of differential forms, one may be of a different opinion with regard to the conservation principle, but this not a point to be discussed here. We just need to realize that Einstein thought in terms of a different mathematical framework, the one of Ricci calculus and covariant and contravariant tensors,  which he had developed in his joint work with M. Grossmann \citep{Einstein/Grossmann}. He started from the observation that a  derivation of the electromagnetic field tensor $F=(F_{ij})$ from a 4-vector potential $\phi=(\phi_i)$ is always possible,  in modernized notation $F=d\phi$. Then Minkowski's equation (\ref{eq Maxwell I Minkoswki}), the first set of the Maxwell equations, boils down to 
\beq d F = dd  \phi = 0 \label{eq Maxwell I Einstein}
\eeq 
 (with $d$ the exterior differential). Of course Einstein wrote this in the form of 4 component equations \citep[equ. (2a)]{Einstein:1916Maxwell}, but still  his equations were easier to grasp and to apply than Minkowski's.
 
For the other set of the Maxwell equations Einstein proposed to ``stay with the generalization of Minkowki's scheme'', i.e. to write it as a divergence equation, although now in  curved spacetime. In order to do so, he introduced the contravariant tensor density\footnote{Einstein called it a ``contravariant $\sqrt{}$-six-vector'' \citep[p. 186]{Einstein:1916Maxwell}.} 
$\mathfrak{F}$ with components
\beq  \mathfrak{F}^{\mu \nu}= \sqrt{|det\, g|}\, \sum_{\alpha, \beta}g^{\mu \alpha}g^{\nu \beta}F_{\alpha \beta} \, .
\eeq 
Then the second set of Maxwell equations became a  generalization of Minkowski's  (\ref{eq Maxwell II Minkowski}),
\beq
\sum_{\nu} \partial_{\nu} \mathfrak{F}^{\mu \nu} = \mathfrak{J}^{\mu}\; , \label{eq Maxwell II Einstein}
\eeq
with $\mathfrak{J}$ the electric current density 4-vector \citep[eq. (5)]{Einstein:1916Maxwell}. Using the covariant derivative $\nabla$ and using the Einstein summation convention, this equation can also be written as
\beq \nabla_{\nu} \mathfrak{F}^{\mu \nu} = \mathfrak{J}^{\mu}  \, . \label{eq Maxwell II Einstein Nabla}
\eeq 
Einstein's proposal was widely influential. One might expect that  Minkowski's dualization idea was filtered out with the transition from the special to the general theory of relativity. 

But this was not the case for all authors. In an important part of the literature on GRT dual complements continued to play a role, although a subordinate one in comparison with other, more central, concepts for the theory. Let us just review a sample of the following three authors, W. Pauli, T. de Donder, and F. Kottler.

{\em Wolfgang Pauli} included a section of its own on  dual complements (``duale Erg\"anzung''),  in the English translation ``dual tensors''  \citep{Pauli:1921English}, in  his well-known review of the theories of relativity for {\em Enzyklop\"adie der Mathematischen Wissenschaften}  \citep[\S 12]{Pauli:1921}. In agreement with his context ($\textsl{dim} M = 4$) he  restricted the discussion to what he called {\em surface} and {\em spatial tensors} (``Fl\"achen- und Raumtensoren''), i.e. to alternating forms of the second and third level $\xi \in \Lambda^k V\; k=2,3$ . He  formulated them already in a differential geometric setting needed for the  general theory of relativity (in later notation $\xi$ in the sections of $\Lambda^2 (TM))$).  Like Sommerfeld, he introduced the duality concept via normalized orthogonal surface elements.\footnote{``With every surface element (\ldots)
in a four-dimensional manifold can be associated another, normal to it,
which has the property that all straight lines in the one are perpendicular
to all straight lines in the other. Such a surface element is called dual to (the first one) if, in addition, it is of the same magnitude.'' \citep[p. 33]{Pauli:1921English}}
 For an alternating 2-tensor $\xi$ with coefficients $\xi^{ij}$ he concluded that the  dual $\xi^{\ast}$  can be given by the corresponding 2-form with coefficients
\beq
\xi^{\ast}_{ij}= \sqrt{|\det \, g|} \,\xi^{kl}\qquad \mbox{with}\quad (ijkl)\sim(1234)\, ,
\eeq
where $\sim$ signifies transformation by an even number of transpositions \citep[eq. (54b)]{Pauli:1921}.  The coefficients of the dual 2-tensor can easily be calculated by ``lifting of indices'' with the metric. 

For an alternating 3-tensor $\xi$, on the other hand,  he similarly arrived at
\beq \xi^{\ast}_i= \sqrt{|\det \, g|} \,\xi^{jkl}\qquad \mbox{with}\quad (ijkl)\sim(1234)\, .
\eeq
Vice versa  the dualization of a vector (field) $\xi$ results in an alternating 3-form
\beq  \xi^{\ast}_{ijk}= \sqrt{|\det \, g|} \,\xi^{l}\qquad \mbox{with}\quad (ijkl)\sim(1234)\, .
\eeq
This allowed to denote integrations of 2- or 3-forms (in later terminology) over surfaces or spatial domains by expressions using tensor or vector densities, i.e. quantities which transform similar to tensors/vectors but with the additional coefficient $\sqrt{|\det \, g|}$ \citep[\S 19]{Pauli:1921}. 

 Dual complements remained for Pauli a subordinate formal tool for specific  calculations. Although he mentioned Minkowski's usage of dual complements in his discussion of the special relativistic Maxwell equations, he  clearly sided with Einstein's preference to do without \citep[\S 28]{Pauli:1921}. In the end, Pauli's adaptation of Grassmann duality to the differential geometric constellation of general relativity remained without further consequences for his overall presentation of the theory. For other authors of the early 1920s, among them     F. \citet{Kottler:1922b} and T.  \citet{deDonder:1923},  this was different. They explored the possibilities for expressing general relativistic Maxwell theory in the terms of dual complements further. 

{\em Th\'eophile de Donder}, e.g.,  adapted Minkowski's presentation of electrodynamics to the context of Einstein gravity . He considered the Maxwell tensor $M$  with three of its components $M_{\alpha\beta}$ encoding the electric displacement and the other three the magnetic field. He introduced the components $M^{\ast}_{\mu \nu}$ (notation by de Donder) of the Faraday tensor as the {\em duals of the functions} $M_{\alpha\beta}$ (``les dualistiques des fonctions'') 
similar to Minkowski. Because of the Lorentzian metric $g_{\alpha\beta}$ of Einstein gravity he explained these ``duals'' by a  formula similar to Pauli's. Rewritten in the light of this commentary his formula boils down to 
\[ M^{\ast}_{\alpha\beta}= (-1)^{\mu +\nu+1}\sqrt{|det\, g|}\,M^{\mu\nu}\, ,
\]
where $\alpha< \beta$ are ``the two numbers of the permutation 1,2,3,4, if one suppresses $\mu, \nu\; (\mu< \nu)$'' \citep[p. 60]{deDonder:1923}.
He then introduced an alternating sum of partial derivatives 
\[ M^{\beta}=(-1)^{\alpha} \partial_{\alpha}M^{\beta\alpha} \;,
\]
which would correspond  to the dual of a Cartanian (exterior) differential $\ast dM$ and  formulated the Maxwell equations in these terms (ibid. p. 64).  In this way de Donder generalized Minkowski's Grassmann duals in the case of 2-forms to the context of  general relativity and applied them in his presentation of electrodynamics  theory, but he did not formulate the Maxwell equations in a form which would boil down to exterior differentials, if expressed in later terminology.

This was different for {\em Friedrich Kottler}. He gave dual complements  a clearer and important conceptual role in a  structural study of  Maxwell's theory. Kottler made a great step towards a foundational study which showed   how far electrodynamics can be formulated without the use of an underlying metrical structure, and at which point metrical aspects come in.\footnote{Today this is called a {\em premetric} approach to electrodynamics. For a recent study deepening Kottler's approach see \citep{Hehl_ea:Kottler}.}
Kottler carved out the dual nature of the two (systems of) Maxwell's equation by opposing the two 4-dimensional field tensors  $F=(F_{ij})$ like in \citep{Minkowski:1908Grundgleichungen} and 
\beq E:=H^{\ast} \, , \label{eq E and Hstar} 
\eeq 
the dual complement of $H=(H_{ij})$ which was his notation for the field excitation, i.e.,  Minkowski's tensor $f$.\footnote{This means $E_{12}=D_3, E_{13}= -D_2, E_{23}=D_1,\; E_{14}=H_1, E_{24}=H_2, E_{34}=H_3$; cf. fn \ref{fn F Minkowski}.} 
All three were (and are) alternating covariant tensors.

 $F$ and $H$ will here be called the {\em first}  and the {\em second}  electromagnetic  field tensors, or also {\em Faraday tensors};   and $E$  the {\em Maxwell tensor}. For $c=1$ the proportionality  (\ref{vacuum proportionality}) in vacuum implies  $H=\mu_0^{-1}F$. For $\mu_0=1$ the two Faraday tensors become (numerically) equal, $F=H$ and accordingly the Maxwell tensor its dual complement, $E=F^{\ast}$.\footnote{Kottler called $F$ ``magnetoelectric six-vector'' and  $E$ as the ``electromagnetic six-vector'' field tensor \citep[pp. 123, 127]{Kottler:1922b}.} 
Because of its role in the source equation (II) of the Maxwell theory,  Kottler's $E$, the Maxwell tensor, is here, like in part of the present literature,  called the  field {\em excitation} and $F$  the field {\em strength}.\footnote{For details see  \citep{Hehl/Obukhov:2003}. }

By playing the game of dual complementation in the inverse direction with respect to Minkowski, Kottler brought the special relativistic Maxwell  equations in a form particularly well adapted to  its invariant  properties with regard to integration. With $S$ the dual complement,  $S= s^{\ast}$, of the ``covariant''  current (i.e. 1-form) $s$ \citep[p. 125]{Kottler:1922b}, 
he  could rewrite Minkowski's variant of the vacuum Maxwell equations  (\ref{eq Maxwell I Minkoswki}), (\ref{eq Maxwell II Minkowski}) with $\mu_0=\epsilon_0=1$ in a form which can be stated  without change of content as
\beq (I) \quad dF = 0 \,, \qquad \qquad (II) \quad dE = dH^{\ast}= S\; , \label{eq Maxwell Kottler}
\eeq 
where $S$ is a current density 3-form. 

Of course Kottler did not use exterior differentials but wrote the equations in terms of sums/differences of partial derivatives \citep[pp. 122, 125]{Kottler:1922b}, e.g. (I) as
\[ \partial_k F_{ml}+ \partial_l F_{mk} + \partial_m F_{kl} = 0 \, ,
\]
 similarly for (II).
But even so he remarked  that the  existence of  (local) 4-potentials $\alpha$ for the  Farady tensor, $F=d\alpha$, and the continuity equation for the 4-current, $dS=0$, are   direct consequences of these equations \citep[p. 127]{Kottler:1922b}. 

As an  advantage of this form   \citet[p. 124f.]{Kottler:1922b} emphasized that the corresponding integral relations  have an  invariant and physically informative meaning. If  we allow us again to use a  modernized notation in terms of integrals of differential forms over a 3-dimensional compact submanifold $A$  of Minkowski space, respectively its  smooth boundary $\partial A$, Kottler could now write the integral version of the  Maxwell equations as:
\beq (I') \quad \int_{\partial A} d F = 0\, , \qquad \qquad (II') \int_{\partial A} H^{\ast} = \int_{\partial A} E  = \int_A dE 
= \int_A S
\eeq
The first electromagnetic  (Faraday) tensor is source free, while the {\em dual complement} of the second  tensor, i.e. the Maxwell tensor for proportionality constants set to 1, is sourced by  the electric  4-current.
 Kottler was convinced that (\ref{eq Maxwell Kottler}) can be considered as the   {\em archetype} (``Urgestalt'') of the Maxwell equations (ibid. p. 129). He took it as  
 a starting point for  exploring  the foundations of electrodynamics also in general relativity. 

He proposed to consider   ``generalized complements'' of   alternating covariant 2-tensors $E=(E_{ij})$ and of alternating 3-tensors $S=(S_{ijk})$ defined with regard to any ``vector of 4-th level''  $e=(e_{1234})$ which ``of course may vary from place to place'',  i.e. an alternating 4-form $e= e_{1234}\,dx^1\wedge\ldots \wedge dx^4$, if we use later notation.\footnote{Kottler used the notation $\epsilon_{1234}$ for $e_{1234}$ which disagrees with the present widespread notation of the ``epsilon-symbol'',  used also below.}
He defined  the {\em generalized complements}   of $E$ and of $S$ as the alternating contravariant 2-tensors $E^{\ast}=(E^{\ast\,ij})$, respectively as contravariant 1-tensor $S^{\ast}=(S^i$) with components given by
\beq E_{ij}= e_{ijkl}E^{\ast\,kl} \; , \qquad \qquad S_{ijk}= e_{ijkl}S^l \, . \label{Kottlers complement}
\eeq
This transformation of an alternating contravariant tensor into an alternating covariant tensor was peculiar for Kottler (neither Grassmann, nor Minkowski or Pauli had considered such a relation).  In the following it will be called {\em Kottler complement}.   

 This feature allowed to formulate an abstract version of the Maxwell equations for generalized Faraday tensors $F$ and $E^{\ast}$ which live in a pre-metric structure characterized by the data $(M,e)$ on a differentiable manifold $M$ and a volume form $e$ \citep[p. 130]{Kottler:1922b}.  Written in  exterior differential notation the  equations  mimicked  (\ref{eq Maxwell Kottler}) closely. (I) remained unchanged,  while (II) had to account for corrections by  a value given in terms of the volume form; and both acquired a generalized meaning:
 \beq  \qquad dF=0\, , \qquad \qquad dE = S \label{eq Kotttlers Maxwell equation}
 \eeq 

After such an extreme generalization of the Maxwell equations he turned towards explaining why  ``a metric enters in the usual presentation'' of electromagnetism: 
\begin{quote}
The reason lies in the   connecting  relations (``Verkn\"upfungsrelationen'') of the two vectors $E$ and $F$.\footnote{``Der Grund hiervon sind die Verkn\"upfungsrelationen zwischen den beiden Vektoren $E$ und $F$.'' \citep[p. 130]{Kottler:1922b}}
\end{quote}
Such a ``connecting relation'' could be linear like in (\ref{eq constants}) or, in general relativity, be established by a metric tensor $g=(g_{ij})$ with volume form $e= \sqrt{|det\, g|}\,dx^1\wedge \ldots \wedge dx^4$.
Then the Kottler complement (\ref{Kottlers complement}) for alternating 2-tensors became  in components, e.g. for $F$ and $F^{\ast}$ \citep[p. 132]{Kottler:1922b},
\beq F_{ij}= \sqrt{|det\, g|} \epsilon_{ijkl}\, F^{\ast kl} \, , \label{eq Kottler complement metric}
\eeq 
with 
\beq
\epsilon_{ijkl} = \Big\{
{   \mbox{\hspace{0.5em}} \mathsf{sg}(ijkl) \quad \mbox{if $(ijkl)$ is a permutation of $(1234)$} 
\atop  \mbox{\hspace{-9em}}  0 \quad \mbox{\hspace{1em}otherwise} }
\label{eq epsilon}
\eeq

In particular for a choice of units such that the numerical value of the permeability of the vacuum becomes $\mu_0= 1$, one can identify  $h=H=F$ (cf. (\ref{vacuum proportionality})).  Kottler's second Faraday tensor then   turned into $E=F^{\ast}$, from which with (\ref{eq Kottler complement metric})
\beq F_{pq}=\sqrt{|det\, g|} \epsilon_{pqkl}\, E^{kl} \, . \label{eq Kottler complement metric B}
\eeq 
Lowering indices, \citet[p. 133]{Kottler:1922b} wrote this relation in the form 
\[ E_{ij}= \sum_{k,l}  \epsilon_{pqkl} (|det\, g|)^{-\frac{1}{2}} g_{ik}g_{jl} \, F_{pq} \, .
 \] 
 In terms of the later Hodge-stgar operator (\ref{eq Kottler complement metric B})  is  nothing but  $F=\ast E$ and thus $E= \pm \ast F$.\footnote{Upper sign for Riemannian, lower sign for  a Lorentzian manifold. }
In other words,  Kottler's complement (\ref{Kottlers complement}) contained  what later would be identified as the Hodge duality for alternating 2-forms or 1-forms in a 4-dimensional (pseudo-) Riemannian manifold as a special case, although he introduced it  in an,  at the outset, premetric approach which presupposed only a (globally given) volume form.

\subsection{\small Cartan \label{subsection Cartan} }
 Kottler was not the only one to realize that a part of relativistic  Maxwell theory could be formulated independent of metrical concepts. {\em \'Elie Cartan} also did so more or less at the same time (1921--1922) when he turned towards studying general relativity from the point of view of differential forms,  infinitesimal group operations and connections. This gave him the opportunity to develop his generalized theory of differential geometric spaces,  which he called ``espaces non-holonomes'' (later Cartan spaces), and a peculiar view of Einstein gravity. He  read Einstein's theory  in analogy with E. and F. Cosserat's generalized theory of elasticity which allowed for hypothetical torque (rotational momenta) in elastic media, in addition to the ordinary stress forces. Motivated by this idea Cartan formulated a generalized view of gravity theory which included torsion in addition to the usual curvature known from Riemannian geometry.  Later this theory  would be called Einstein-Cartan gravity.\footnote{For a survey of Einstein-Cartan theory see \citep{Trautman:EC}, for the development of Cartan's thought  
\citep{Chorlay:Cartan,Nabonnand:2016Cartan,Scholz:2019Cartan}. \label{fn Chorlay et al}}

Cartan knew Grassmann's work well,  he appreciated in particular Grassmann's combinatorial products. They  helped him to formulate his calculus of  differential forms, e.g. $\omega= \sum_i a_i dx_i, \, \tilde{\omega}= \sum_i b_i dx_i$,   and in particular the  alternating product of differential forms,\footnote{Cf. \citep{Katz:James,Katz:1985}.}
 which in his notation appeared as
\[  [\omega \, \tilde{\omega}] = \sum_{i<j} (a_i b_j - a_jb_i)[dx_i\,dx_j]
\]
He also introduced the exterior differential $\omega'$ of a form $\omega$,  e.g.,
\[ \omega'= (\frac{\partial a_i}{\partial x_j}-\frac{\partial a_j}{\partial x_i})\, [dx_i\, dx_j] \, .
\]
In the following  Cartan's differentials $\omega'$ will  be written in the more recent form $d\omega$ in place of  $\omega'$, while we stick  to his  () notation of the alternating product by square bracket without comma.

Cartan's  described the  infinitesimal neighbourhoods (later understood as tangent spaces $T_pM$)  of a point $p$   with coordinates $x_1, \ldots, x_n$ by using differential forms  $\omega_1, \ldots, \omega_n$ (later understood as the  elements of a coframe of $T_pM$) and  characterized the geometry in the neighbourhood by infinitesimal group operations which mimicked, in a specific way, the Kleinian approach to geometry with a principal group $G$. This led him to characterize the geometry by a {\em connection} in the group, given by a set of differential 1-forms  $\omega_{kl}$. The indices $k,l$ indicate the components of the infinitesimal group. As as a whole, i.e. considering the complete matrix $(\omega_{kl})$, this collection can be understood as a 1-form with values in the infinitesimal group (the Lie algebra) of $G$ which served as a  generalized rotational group for the geometry. Moreover he also  considered the coframe given by the $\omega_j$  as  a kind of translational connection, complementary to the rotational connection given by the $\omega_{kl}$. He analysed the deviation of such a structure from the corresponding (global) Kleinian geometry and characterized it by two types of curvatures, the rotational curvature $\Omega_{kl}$ and the translational curvature $\Omega_i$. Both were constructed as differential 2-forms, the second one as a whole  vector valued $(\Omega_i)$, and the first one, considered as a matrix $(\Omega_{kl})$ with values in the Lie algebra of $G$. 
This opened up the way to studying a new class of differential geometric spaces, later called Cartan spaces.\footnote{For a recent mathematical textbook see \citep{Sharpe:Cartan_Spaces}, for  historical literatur  see fn. \ref{fn Chorlay et al}. } 
In the papers considered here he considered special orthogonal groups of Euclidean or Lorentzian signature in the dimensions $n=3, 4$..

Cartan made use of Grassmannian complements at different places in his work, most clearly in his large study {\em Sur les vari\'et\'es \`a connexion affine et la th\'eorie de la relativit\'e g\'en\'eralis\'ee} \citep{Cartan:1923/24,Cartan:1925}.\footnote{Reprint of both parts in \citep{Cartan1955}, English in \citep{Cartan:1986}.}
He  introduced the dualization of  alternating products in an exemplary way in a  case by case discussion. 

In the 3-dimensional case, for example, described by mobile orthogonal but not normalized 3-frames $(e_1,e_2,e_3)$ and the associated  metric $g = diag\, (g_{11}, g_{22}, g_{33})$, he mentioned  in passing that any {\em bivector} (a Grassmannian quantity of the second level) can be represented in the form
\[  \Omega^{13}\,[e_1\, e_3] +  \Omega^{31}\,[e_3\, e_1] + \Omega^{12}\,[e_1\, e_2]\, ,\]
or just as well by {\em polar vector} of the same measure (``de la m\`eme mesure'') \citep[\S 79, p. 16]{Cartan:1923/24}
\[ \frac{1}{\sqrt{g_{11}g_{22}g_{33}}}(\Omega_{23}\,e_1 + \Omega_{31}\,e_2+ \Omega_{12}\,e_3 ) \, .
\]
This generalized the Grassmannian view of the relation between ``axial'' and ``polar'' vectors to differentiable manifolds. 

Clearly this kind of relationship was reciprocal (ibid, \S 60, p.400). Similarly, in  dimension $n=4$ there exists an ``invariant correspondence'' between a vector (field) $\xi^i e_i$ and a ``polar trivector''
\[ \sum_{i=1}^4 \epsilon_{(ijkl)} \, \xi_i\, [e_j\, e_k\, e_l]\qquad \mbox{(where $j<k<l$)} \, ,
\]
where the symbol $\epsilon_{(ijkl)}$ is used like in (\ref{eq epsilon}); similarly  also between ``a system of bivectors'' $\xi^{ij}\, [e_i\, e_j]$ and the system of ``polar bivectors''\footnote{Cartan wrote the sum term by term with explicit signs \citep[\S 60, p. 401]{Cartan:1923/24}.}
 \[ \sum_{i<j} \epsilon_{(ijkl)}\, \xi_{ij}\, [e_k\, e_l] \qquad \mbox{(with $k<l$)} \, .
 \]
We have seen above that in other parts of the literature  such Grassmann-type reciprocal relationships were called {\em dual complements}. Although it was not Cartan's  terminology we will call them so in the sequel.   
 For Cartan's interpretation of the Einstein tensor, for his generalizations of the theory of gravity, and in the discussion of the relativistic Maxwell equation dual complements turned out to be of central importance.

When Cartan started studying  Einstein's theory he looked for an analogy to the generalized theory of elasticity  proposed by the brothers E. and F. Cosserat.\footnote{See \citep{Brocato/Chatzis:Cosserat,Scholz:2019Cartan}.}
 Cartan developed a peculiar interpretation of the Einstein equation which in ordinary symbolism of differential geometry is
\beq G = \kappa\, T \; ,\label{Einstein eq}
\eeq
with the Einstein tensor $G=Ric-\frac{R}{2}g$ formed from curvature expressions  of spacetime, the stress-energy-momentum tensor  $T$, and the gravitational constant $\kappa$  of Einstein theory.\footnote{$Ric$ denotes the Ricci tensor $g$ the metric, $R$ the scalar curvature of spacetime, $\kappa= 8\pi G_N$ with the Newton constant $G_N$ and $c=1$.}
Cartan considered the right hand side of (\ref{Einstein eq})  as a quantity describing the dynamical state of matter (``quantit\'e mouvement mass'' \citep[\S 78, p. 13]{Cartan:1923/24}). The equation itself seemed  to allow him  an identification of a geometrical curvature quantity with a matter quantity.\footnote{``On sait que, dans la th\'eorie de la relativit\'e g\'en\'eralis\'ee d'Einstein, le tenseur qui caract\'erise compl\`etement l'\'etat de la mati\`ere au voisinage d'un point d'Univers est {\em identifi\'e} \`a un tenseur faisant intervenir uniquement les propri\'et\'es {\em g\'eom\'etriques} de l'Univers au voisinage de ce point'' \citep[p. 437, first emphasis ES, second emphasis in the original]{Cartan:1922[57]}.}

As already indicated,  in  Cartan's approach the curvature  was expressed by a set of differential 2-forms $\Omega_{ij}$ the indices of which indicate components of  the ``infinitesimal rotations'' (the Lie algebra of the special orthogonal group).\footnote{For an analogue to the Einstein equation in  dimension $n=3$ he considered the rotation group of Euclidean space, for $n=4$ the one of Minkowski space, i.e. the Lorentz group.}
He rewrote the Einstein tensor  as an $(n-1)$-form in his curvatures (rotational and translational) for $n=3$ and  $n=4$. Here we need not discuss  the derivation of the curvature  forms and Cartan's way of writing the Einstein tensor.\footnote{For more details see \citep{Scholz:2019Cartan}. Be aware of the error note for the journal version of this article. The arXiv version is correct.}
The stress-energy tensor $T=(T^i_{\;j})$ of Einstein theory is  a vector valued 1-form (for both $n=3,4$). A transformation from  $(n-1)$-forms to 
1-forms was possible, in principle, by means of a Grassmann type duality. 
Cartan, however, did not care about such a move; he found it more natural to express the dynamical state of matter by a 2-form ($n=3$) or 3-form ($n=4$) respectively, which assigns forces and rotational momenta to a surface or spatial element of spacetime. 
But the infinitesimal rotations and translations (the values of the curvature forms $(\Omega_{kl})$ and $(\Omega_i)$) had to be transmuted into vectors or momenta, which could be interpreted as the stress force or torque acting on the respective surface element. 

For  dimension $n=3$ this was easy to achieve. Infinitesimal rotations in dimension $n=3$ can be represented in vector analysis by vector products anyhow; they could  be associated to a vector (at least implicitly by a Grassmann duality as we have seen in section \ref{subsection vector operations}). On the other hand, the infinitesimal translations of $(\Omega_i)$ could be transmuted into a ``bivector'', a Grassmannian quantity of the second level. Then they could be interpreted as the geometrical expression of  rotational momenta. This was the background for Cartan's, at first sight surprising,   choice to call the {\em  translational curvature} of his new type of geometry  {\em torsion}, which he announced already in his {\em Comptes Rendus} notes:
\begin{quote}
To any closed infinitesimal loop there is generally an associated translation; in this case  one can say that the given space  differs from Euclidean space in two respects: 1. by a \underline{curvature} in the sense of Riemann, which is expressed by the rotation; 2. by a \underline{torsion} which is expressed by the translation (emphasis in the original). \citep[594f.]{Cartan:1922[58]}\footnote{``Dans les cas g\'en\'eral o\`u il y a une translation associ\'ee \`a tout contour ferm\'e infiniment petit, on peut dire que l'espace donn\'e se diff\'erencie de l'espace euclidien de deux mani\`eres: 1$^{\circ}$ par une {\em courbure}  au sens de Riemann, qui se traduit par la rotation; $2^{\circ}$ par une {\em torsion}, qui se traduit par la translation'' \citep[594f.]{Cartan:1922[58]}.
}\end{quote}

 For $n=4$ the case was more complicated.  In a long  discussion  Cartan found that there arose a problem for the enhancement of  the  curvature by a translational component \citep[\S 78--\S 83]{Cartan:1923/24}. In fact, it would be inconsistent with the Maxwell-Lorentz theory of electrodynamics if one did not assume an additional term for the  mass-energy (``quantit\'e mouvement mass''). 
  This led him to adding an expression   to the the energy momentum tensor (respectively 3-form),   which much later would become known as a spin term.\footnote{See the commentary by A. Trautman in \citep{Cartan:1986} and with regard to history \citep{Scholz:2019Cartan}.}

We need not go here to the heart of the problem Cartan found in the Maxwell-Lorentz theory for his generalized theory of gravity. For our purpose it will suffice to shed a glance at his discussion of the relativistic Maxwell equation. 
To start with,  he rewrote special relativistic Maxwell theory coherently in terms of differential forms and exterior derivatives. In particular, the (first) Faraday tensor and the Maxwell tensor were represented  as 2-forms $\Omega, \overline{\Omega}$ with 
\beqarr \Omega &=&  \sum_{(ijk)}B_i\, [dx_j\,  dx_k] +  \sum_i E_i\, [dx_i\, dt] \\
 \overline{\Omega} &=& \sum_{(ijk)}D_i\, [dx_j\,  dx_k] +  \sum_i H_i\, [dx_i\, dt] \, ,
\eeqarr
where $(ijk)$ indicates the summation over all cyclic permutations of $(123)$, and the charge/current density $S$ as a 3-form
\[ S = \rho\, [dx_1\, dx_2\, dx_3] - \sum_{(ijk)} J_i\, [dx_j\, dx_k\, dt]
\]
 Taking account of his sign conventions the  Maxwell equations came out as
 \beq d\Omega = 0 \, , \qquad \qquad d\overline{\Omega}=- 4\pi S \; , \label{eq Cartan Maxwell}
 \eeq 
 that is like Kottler's equation (\ref{eq Kotttlers Maxwell equation}), up to sign and with units sucht that numerically  $\mu_0 = 4\pi$ \citep[\S 80]{Cartan:1923/24}.
 
 Because of his method of orthonormal frames it was not too difficult for Cartan to import  Maxwell theory into a general relativistic framework. Assuming $n=4$ differential forms  $\omega^0, \ldots, \omega^n$ which at every point $p$ represented an orthonormal Lorentzian  co-frame to some basis (frame) $e_0,\ldots, e_3$ of the infinitesimal vector space at $p$.\footnote{Cartan spoke of ``syst\`eme de r\'ef\'erence de Galil\'ee'' \citep[\S 82]{Cartan:1923/24}.}
 the Faraday and Maxwell tensors of special relativity turned into 2-forms on the generalized space, which  in slightly adapted notation were 
\beq  \Omega = \sum_{i<j} H_{ij}\, [\omega^i\, \omega^j]\,, \qquad \qquad \overline{\Omega} = \sum_{i<j,\, k<l} \epsilon_{ijkl}\, H^{ij}\, [\omega^k\, \omega^l] \, ,
\eeq
where $i,j,k,l \in {0, \ldots 3}$ and $\epsilon_{ijkl}$ like in (\ref{eq epsilon}).
 Apparently Cartan  understood the lifted  indices for $H_{ij}$ like in  the Ricci calculus. In  his orthonormal frames the metric is diagonalized and of  Lorentzian  signature,  $g= diag\, (+1,-1,-1,-1)$. The lifting of indices thus has only consequences for the sign of the coefficient. 
 $\overline{\Omega}$ denoted thus a Grassmann type dual complement  written in Cartan symbolism. He concluded that with an adapted definition of the electric current,
 \[ S = \sum_{ijkl} \epsilon_{ijkl} J^i\, [\omega^j\, \omega^k \, \omega^l] \, ,
 \]
 the Maxwell equations  (\ref{eq Cartan Maxwell}) are still valid.  Moreover 
 \begin{quote}
 \ldots  they don't make use of the affine connection of the universe \citep[\S 82, p. 20]{Cartan:1923/24}.
 \end{quote}
So Cartan was well aware of the fact that  Maxwellian electrodynamics depends only weakly on the metric of spacetime, and not at all on the affine connection. The equations themselves can be formulated in terms of exterior differentials; they do not depend on the metric at all. The  interrelation between the Faraday tensor and the Maxwell tensor can be expressed in terms of a Grassmann type duality, later called Hodge duality. It thus depends indirectly on  metrical concepts of spacetime, although for the case $n=4$ only on the conformal class of the metric.   Moreover we have  seen, how Grassmann type complements were used by  Cartan  also for other purposes in his investigations of gravity.

\section{\small The birth of Hodge duality \label{section Hodge duality}}
It is now time to turn towards the second main topic of our historical report, Hodge theory.  
{\em William V.D. Hodge} (1903--1975) started the study of algebraic surfaces and varieties at the end of the 1920, at a time when modern methods of algebra were increasingly introduced into algebraic geometry and the combinatorial topology  ({\em analysis situs}) of manifolds and complexes. He became a central figure for enriching this field by generalizing Riemann's transcendent methods  in the study of compact Riemann surfaces 
to higher dimensions. Before we look at his work we have to shed a short glance at the background knowledge from which he started.

\subsection{\small Background: Riemann's  theory  of Abelian integrals  \label{Background Riemann}}
 Riemann had taken advantage of the fact that the real and imaginary parts $u\,dz$ and $v\,dz$ of a holomorphic (or meromorphic) form $\omega = u dz + v\,\mathrm{i} \,  dz$ ($z$ complex coordinate)  are harmonic, i.e., the coefficient functions satisfy the partial differential equation $\Delta u=0, \,\Delta v = 0$. In his influential paper on the theory of Abelian functions \citep{Riemann:1857} he characterized the (real) harmonic  forms on a  Riemann surface $\Fs$, after cutting the latter along a number of curves into a simply connected surface, by  boundary value problems and proposed to solve them by applying the Dirichlet principle.  At his time this method was not yet mathematically well developed, and it took half a century and the work of many mathematicians to fix the loose ends. At the end point of this story stands Hilbert's  vindication of the Dirichlet principle  \citep{Hilbert:1904} and Weyl's famous book on the {\em Idee der Riemannschen Fl\"ache} \citep{Weyl:IdeeRF}.\footnote{The problems of this approach and its solution have been discussed historically at many places, see in particular \citep{Monna:Dirichlet},  \citep[chaps. 5, 6, 7.7]{Bottazzini/Gray}, \citep[chaps. 16--18]{Gray:Analysis}.}
A crucial outcome of his approach were deep  insights into the interrelation of the topology of a compact Riemann surface $\Fs $, its complex structure, and its algebraic description. In particular Riemann studied closed curves (by later authors called ``retrosections'')  which  neither partially nor in  total bound a  part of the surface and found that their maximal number is even, $2p$
   \citep[\S 3, p. 104]{Riemann:1857}. He  called   $2p+1$ the ``order of connectivity'' of $\Fs $. Later one would speak of $b_1=2p$  as the first Betti number of the surface.

He studied  integrals  $w_j= \int_{\gamma}\omega_j$ over  a curve $\gamma$ (closed or not) with holomorphic $\omega$, so-called  {\em Abelian integral of the first kind}   (of the second and third kind for meromorphic differentials without, respectively with poles of first order), and found that there are only finitely many  linearly independent ones $w_1, \ldots, w_q$ on a given surface $\Fs$  \citep[\S 4, p. 105]{Riemann:1857}.\footnote{The corresponding  differential forms  $\omega_1, \ldots, \omega_q$ are  independent if no linear combination of them is a total differential $dW$ of a univalued function $W$ on $\Fs$.}  
Due to his usage of the Dirichlet principle for the real and imaginary parts of Abelian integrals  $w= u +  \mathrm{i} \,  v$, the maximal number of independent holomorphic forms turned out to be determined by the structure of dissection of $\Fs$ into a simply connected surface. It agreed   with the number identified in the analysis situs study of the surface, $q=p$. Riemann considered such an integral $w_j$ as a  multivalued function on $\Fs$,  the values of which differ by (integer) linear combinations of  $a_{ij}=\int_{\gamma_i}\omega_j\;, i=1,\ldots , 2p$, the ``Periodicit\"atsmoduln'' \citet[p. 105 etc.]{Riemann:1857} later simply called {\em periods}.

 He also showed that a  compact Riemann surface can be represented as a  complex  {\em algebraic  curve} with  an algebraic equation $F(z,w)=0$ in two complex variables $z, w$  in various ways. He argued that two representations of this type can be transformed  into another algebraically (by what later would be called birational transformations) and  showed that for  a curve of order $n$ with $w$ simple branch  points (and no ones of higher order) the number $p$ of Abelian integrals satisfies the relation
 \beq p=\frac{w}{2}- (n-1) \label{eq genus} 
 \eeq 
 and can thus be read off from the algebraic representation of the Riemann surface as a curve  \citep[\S 7, p. 114]{Riemann:1857}.\footnote{See also the annotation (2) by H. Weber in 1876, the editor of the first edition of Riemann's {\em Werke}.}

This opened the way for an algebraic study of compact Riemann surfaces,  viz. complex algebraic curves embedded in the complex projective plane. It was pursued  further by A. Clebsch, M. Noether,  A. Brill and many other mathematicians. These authors were able to replace Riemann's analytic (``transcendent'') methods by algebraic ones and characterized in particular the Abelian differentials/integrals by algebraic expressions (``adjoint polynomials/curves''). In this context  Clebsch  introduced the name {\em genus} (``Geschlecht'') for the algebraically determined  number (\ref{eq genus}). For him and mathematicians in his circle  it must have been clear (already in the 1860s and 1870s)  that this number was a birational invariant and that it characterized the analysis situs (topological) property of the corresponding Riemann surface.
F.Klein explained this connection between different aspects of Riemann's theory  to a broader (mathematical) audience \citep{Klein:RiemannscheFlaechen}.\footnote{\citet{Le:genre} also emphasizes this role of Klein's book, but imputes that for Clebsch and the mathematicians of his generation the interrelationship between the different aspects of Riemann's theory was unclear or doubtful. I agree with the first characterization, but not with the second.}

The situation became much more complicated, when  attempts were being  made to apply Riemann's integrated approach to  the study of complex algebraic surfaces.
\citet{Clebsch:1868} introduced the {\em genus} $p$  of an algebraic   surface $\Ss$ possessing  only simple singularities as the number of linearly independent double integrals of a special type, defined  algebraically.  \citet{Cayley:1871} attempted to find a formula analogous to (\ref{eq genus}) for the genus of a surface, but realized that the number so defined did not always agree with Clebsch's value $p$. M. Noether therefore called it the {\em arithmetical genus} $p_a$ in contrast to Clebsch's, which he now called {\em geometrical genus} $p_g$. Both types of genus turned out to be birational invariants.\footnote{See \citep[p. 313f.]{Brigaglia_ea:Abel}.} 

About this time \citet{Betti:1871}  indicated how Riemann's analysis situs concepts, in particular the order of connectivity, can be generalized to higher dimensions. For a manifold  of dimension $n$ he introduced connectivity numbers $b_k$ (later called Betti numbers) for any  dimension $1\leq k\leq n$. But it remained completely unclear whether or how the different types of genera of a complex algebraic surface $\Ss$ (with real dimension $dim_{\R}\, \Ss = 4$) had something to do with  these topological invariants (the $b_k$). 
From now on the different aspects of Riemann's integrated research program combining analysis situs, algebraic geometry, differential forms and their integrals (even specified to the case of holomorphic or meromorphic form), which to a certain degree even 
tied up with   the differential geometry of manifolds,\footnote{Cf. \citep{Scholz:Diss}.} 
started to be  developed  in different research traditions with weak overlap or interchange. This remained  essentially so  for about half a century. 

\subsection{\small Hodge's first steps towards generalizing  Riemann's theory  of Abelian integrals \label{Hodge first steps}}
When Hodge entered research in  mathematics, several of the fields treated by Riemann  had been developed further and stood at the brink of becoming mathematical subdisciplines of their own. Hodge assimilated a wide area  of literature with different approaches to what would become the geometry and topology of the 20th century.\footnote{See, e.g., the literature cited in \citep{Hodge:1932Dirichlet,Hodge:1932/34}.}
Riemannian differential geometry and its generalizations  got an immense push with the rise of the general theory of relativity.\footnote{See among others \citep{Bourgignon:connexion,Bottazzini:Connection,Scholz:Connections,Scholz:2019Cartan,Reich:Tensoren,Reich:Connection}}
 On the other hand,  the algebraic geometry of complex varieties had accumulated a rich corpus of insights in particular established by the Italian school of geometers; although  it was not always based on reliable foundations.\footnote{\citep{Brigaglia/Ciliberto,Brigaglia_ea:Abel,Schappacher:2015}}

The analysis situs study of manifolds and complexes, at the  beginning of the new century often called combinatorial topology, was being reshaped by modern algebraic concepts and turned into   {\em algebraic topology}.\footnote{See  \citep{James:History,James:alg_top,Epple:Knoten,Herreman:1997,Herreman:2000,Scholz:James,Volkert:Homeomorphieproblem} and the respective contributions to this volume.} 
The integration of differential forms over submanifolds, in particular closed ones (so-called ``cycles'') had been introduced in higher dimensions by E. Picard, H. Poincar\'e  and E. Cartan. By generalizing the theorem of Stokes these authors realized that the integral  $\int_{c} \omega$ of a  differential form $\omega$ of degree  $k$ with vanishing  derivative ($d\omega=0$), taken  over a closed submanifold $c$ with $dim\, c=k$ 
is the same for any homologically equivalent submanifold $\tilde{c}$ , $\int_{\tilde{c}} \omega= \int_{c} \omega$. Thus the  {\em periods} $\int_{c_j} \omega$ of such a  differential form  with respect to a generating system $c_1, \ldots, c_m$ of the $k$-th homology could  be considered as belonging to the analyis situs of the manifold, long before the idea of  cohomology theory was shaped. During the 1920s E. Severi and S. Lefschetz were the main protagonists of using this type  of analysis situs for the study of algebraic surfaces.

At the turn to the 1930s \citet{deRham:1931} made in his PhD dissertation   a decisive step forward towards establishing a  (dualizing) analogy between forming the  boundary $\partial \, c$   of a (differentiable) complex $c$  and  the exterior derivative $d\omega$ of a $k$-form $\omega$.\footnote{For de Rham and his environment see \citep{Chaterji_ea:2013}.}
 He  introduced homological terminology into the treatment of differential forms: $\omega$ is {\em closed} (``ferm\'ee'') if $d\omega=0$ (written by him in Cartan's notation of the time $\omega'=0$); it is {\em homologue zero}, $\omega \sim 0$  (notation used by de Rham),  if it is the exterior differential of another form $\widetilde{\omega}$ of degree $k-1$ , $\omega=  \widetilde{\omega}'$ etc. (p. 176). He showed: 

{\footnotesize
\begin{itemize}
\item[--] Every closed form is homologue to a linear combination of finitely many ``formes \'el\'ementaires'' \citep[p. 180]{deRham:1931}.
\item[--]  A closed form with all periods zero is itself homologue zero,  $\omega\sim 0$  (p. 185).
\item[--] Given  $m$ homologously independent $k$-cycles $c_j$ and $m$ rel values $a_1, \ldots , a_m$ on can find a closed $k$-form $\omega $ such that $\int_{c_j}\omega = a_j$ for $1\leq j \leq m$ (p.186).
\item[--] And finally, the maximal number of homologously independent closed $k$-forms $q_k$ coincides with the $k$-th Betti number $q_k=p_k$ (p. 187).
\end{itemize} }
 In this sense  de Rahm established the basic insights into what later would become the (de Rham) cohomology of differentiable manifolds, although he stopped short of introducing the cohomology groups of differential forms themselves.\footnote{See \citep{Katz:1985,Katz:James,Massey:cohomology}.}

In the late 1920s Hodge  started  studying  the periods of rational forms on algebraic varieties. He  joined and expanded the research program of E. Severi and S. Lefschetz, who promoted a (re-) integration of analysis situs methods into the theory  of complex algebraic surfaces and also of algebraic varieties $V$ of higher dimensions. In these attempts 
 integrals  $\int_{c} \omega$ like above were studied, but here with rational $k$-forms $\omega = R(x_1,\ldots, ,x_k)dx_1 \ldots dx_k$ which lead to finite integrals over analytically defined subvarieties or even complexes $c$ (``complexes analytiques'') of  dimension $dim_{\R}\, c = k$. Their integrals, called of the {\em first kind} like in the case of algebraic curves,   were expected to give important information on the variety $V$; but their study led to 
much more complications   than in Riemann's context.\footnote{``On ne sait pas grand' chose sur ces int\'egrales, mais il est probables qu'elles ont une importance consid\'erable pour la th\'eorie alg\'ebro-arithmetique de $V_d$'' \citep[p. 55]{Lefschetz:1929}. Lefschetz considered the number $p_g^{(k)}$ of independent $k$-fold integrals of the first kind as  a generalization of Clebsch/Noether's geometric genus $p_g$ of curves.}
Already for  $n=2, \, k=2$ it was unclear whether there are (non-vanishing)
integrals of the first kind, all  periods of which are zero. Severi asked for their number and neither he nor Lefschetz expected  that it might be zero by general reasons.\footnote{\citep[p. 350]{Lefschetz:1921}, \citep[p. 28ff.]{Lefschetz:1929}.} 

Using  analytic considerations, i.e. by using ``transcendent'' methods  rather than relying on exclusively algebraic ones, Hodge concluded   that this number is zero; but with the present analytic tools at hand  a proof was difficult \citep{Hodge:1930}. Only after an initially  strong opposition  Lefschetz accepted Hodge's claim \citep[p. 175f.]{Atiyah:1976Hodge}. For Hodge the dissertation of de Rham came just at  the right time for 
 evolving his approach. From 1931 onward he could build upon the  methods of general  differential forms and concentrate on  specifying them for the case of {\em harmonic forms}, in order to   generalize Riemann's analytic theory to higher dimensions.  But how could  one characterize  harmonic  forms for complex algebraic varieties of dimension $n>1$? If one wanted to use the  Beltrami-Laplace operator  $\Delta$ for characterizing harmonic functions $u$,
\beq \Delta\, u = \frac{1}{\sqrt{ det\, g}} \, \partial_i (\sqrt{ det\, g} g^{ij}) \partial_j u = 0 \; , \label{eq Beltrami-Laplace op}
\eeq 
 a Riemannian metric 
$g=(g_{ij})$ on the manifold  had to be presupposed.\footnote{In modernized notation the formula boils down to $\Delta u = \nabla_i \partial^i u =0$  with  $\nabla$  the covariant derivative  associated to $g$.} 

In his first sketches of his theory and the announcement of it  \citet{Hodge:1932Dirichlet} avoided the problem and argued with  ``euclidean $n$-cells'', i.e.  he assumed the cell complexes studies as embedded in a Euclidean space. But already in the paper \citep{Hodge:1932/34} written in 1932, although published only two years later, and in an outline of his theory \citep{Hodge:1933MathGaz} he explained how a complex algebraic variety $V$ of dimension $m$ could be endowed with a Riemannian metric. He assumed that $V$ can be given in a singularity free form and  embedded in a projective space of sufficiently high dimension $r$.\footnote{``The general theorem has not yet been 
proved, but we shall assume its truth, or better, we shall confine 
our attention to varieties which can be transformed into varieties
without multiple points, and we shall suppose that $V$ is a variety
of $m$ dimensions without singularities, lying in a complex projective
space of $r$ dimensions $(z_0, \ldots z_r)$'' \citep[p. 304]{Hodge:1933MathGaz}.}
He then used a method by G. Mannoury for endowing the projective space $P^r(\C)$, and with it also the subvariety $V$, with a Riemannian metric.\footnote{\citep{Mannoury:1900}}
 Mannoury proposed to  embed  the projective space  $P^r(\C)$  in an  Euclidean space of even higher dimension $s=(r+1)^2$ and to use the induced metric of the embedding $P^r(\C) \hookrightarrow \E^s$ in the Euclidean space of dimension $s$.\footnote{With $X_h,X_{hk}=X_{kh}, Y_{hk}=-Y_{kh}, \; (h,k=0, \ldots r)$ coordinates of $\E^s$,  one sets $X_h= \sqrt{2 z_h\overline{z}_h}$, $X_{hk}=z_h\overline{z}_k + z_k\overline{z}_h, Y_{hk}= \mathrm{i} (z_h\overline{z}_k - z_k\overline{z}_h)$, where the projective coordinates are constrained by $\sum_{j=0}^r |z_j|^2=1$ \citep[p. 304]{Hodge:1933MathGaz}.   }

 Hodge  apparently found it clear that different, but birationally equivalent presentations of $V$  ought to  lead to the same invariants for the algebraic manifold. At least he did not hesitate to speak of the submanifold $M$ which corresponds to the points of $V$ as the {\em Riemannian manifold} of $V$ and considered $M$ as a higher dimensional analogue of the Riemann surface $\Fs$  of an algebraic curve. In this expectation Hodge turned out to be right, but it took a long way to go until this conjecture could be justified. 
In his look back Atiyah emphasizes the problem by talking  of Hodge's 
``apparently strange idea of introducing an auxiliary metric into algebraic geometry'', which could be vindicated only much later.\footnote{ According to  \citet[p. 187f.]{Atiyah:1976Hodge} a  first vindication of Hodge's application of his theory  to algebraic manifolds resulted from   his  proof that the  decomposition of the space of harmonic forms $\Hs^r = \sum_{p+q=r}\Hs^{p,q}$ (see below) the dimensions $h^{p,q}= \mathrm{dim}\, \Hs^{p,q}$ are invariants of the complex structure of $V$. An intrinsic definition of the  $h^{p,q}$, the Hodge numbers, became available only in the 1950s after the introduction of sheaf theory (see sec. \ref{Hodge theory in sheaf cohomology}). }
 The problem is:
 \begin{quote}
 For Riemann surfaces the complex structure defines a conformal structure and hence the Riemannian metric is not far away, but in higher dimensions this relation with conformal structures breaks down and makes Hodge's success all the more surprising.  \citep[p. 187f.]{Atiyah:1976Hodge}
 \end{quote}
It thus was a daring move, combined  perhaps with  a  visionary perspective and a bit of luck, which  allowed Hodge to  anticipate crucial insights into the role which harmonic differential forms and their integrals could play as an intermediary between topology and complex  algebraic geometry.

\subsection{\small The $\ast$-operation, the Hodge theorem and  algebraic surfaces \label{subsection Hodge star-operation}}
Early in the 1930s our author announced a central  theorem of his new theory and showed how it could applied to the study of algebraic surfaces \citep{Hodge:1933genus}.
\begin{theorem}
On an orientable Riemannian manifold \footnote{In Hodge's formulation, an ``analytic construct of $n$ dimensions which has the topological properties of an orientable absolute manifold \ldots which has attached to it a Riemannian (positive definite metric)'' \citep[p. 312]{Hodge:1933genus}.}
with  Betti numbers $p_m$ there are exactly $p_m$ harmonic $p$-forms $ (1\leq m \leq n)$. 
\end{theorem}
This was the first version of what later would become known as Hodge's theorem. 
The proof followed in two papers in the {\em Proceedings of the London Mathematical Society}  \citep{Hodge:1932/34,Hodge:1934Harmonic}. In these papers Hodge mentioned de Rham's dissertation\footnote{\citep[p. 257]{Hodge:1932/34}  \citep[p. 90f.]{Hodge:1934Harmonic}}  and profited from its insight with regard to the dual character of (closed) differential forms and the corresponding cycles in the manifold. He introduced a  basis $\{w_i\}$  of the harmonic forms of degree $m$, related  to  a basis of the $m$-th homology represented by cycles $\{\Gamma_j\}$ ($  1\leq i, j \leq p_k $),  by the condition $\int_{\Gamma_j}w_i=\delta_{ij}$. If we use  the later notation $H_k(M,\R)$ for the homology and $\Hs^m(M)$ for the harmonic forms (of the first kind) of degree $k$ on   the manifold $M$, Hodge  clearly noticed a duality relation between the two, even though he did not yet use  the word.
In a follow up paper read in February 1934  he took stock of what was achieved    and set out to give ``an account of the  principles on which the method is based''   \citep[p. 249]{Hodge:1935HarmonicIntegrals}.

If we use an ex-post notation of de Rham's insight into the relation between ``homologue differentials'' and ``cycles'' as a duality relation between  homology  and harmonic forms
we can resume {\em Hodge's theorem} in retrospect  as
\beq \Hs^k(M) \cong H_k(M,\R)^{\ast}  \; . \label{eq Hodge's thm}
\eeq

In the first paper, \citep{Hodge:1933genus},   the general theorem was stated without using the terminology of {\em harmonic forms}, while it appeared prominently in the  next publications. 
In 1933 Hodge rather used the language of   ``linear independent skew symmetric tensors'' $B_{i_1\ldots i_m}$  satisfying two equations (1), (2).
Using an abbreviated notation $B$ for Hodge's  alternating tensor (or form), his  equation (1)  expressed a vanishing  exterior derivative, $dB=0$. The second one demanded  the vanishing of ``contravariant derivative'', i.e. a covariant derivative with lifted index $\nabla^r B = 0$. Only in the review paper    \citep{Hodge:1935HarmonicIntegrals} he introduced an equivalent formulation of equation (2) by expressing it 
 as the vanishing of the exterior derivative of the Hodge dual $d\ast B=0$ (see below).\footnote{
The tensor equations in \citep[p. 312]{Hodge:1933genus} (with original equ. numbers)  were:
\beqa \sum_{r=1}^{m+1}(-1)^{r-1}B_{i_1\ldots i_{r-1} i_{r+1} \ldots i_{m+1}, i_r} &=& 0 \qquad (1) \\
g^{rs}B_{i_1\ldots i_m, s} &=& 0 \qquad (2)
\eeqa
}
In the 1933 paper he introduced  the terminology ``harmonic integrals'' and ``harmonic forms''  only a  in the context of  double integrals on algebraic surfaces ($n=4, m=2$), \citep[p. 315ff.]{Hodge:1933genus}.

For a detailed  proof of the announced theorem he referred to  \citep{Hodge:1932/34}; in 1933  he gave only an  outline of the 
argument, but  emphasized   its importance for understanding  the geometric genus $p_g$ of a (singularity free) algebraic surfaces $\Ss$. 
We cannot go  into details of this interesting history; but  certain aspects of Hodge's discussion of algebraic surfaces shed light on our story. Here, for dimension $n=4$ and $m=2$, he showed that for a skew symmetric  tensor  now called  $\phi= B_{ij}dx^i dx^j$ satisfying the two equations (1), (2) of his theorem another form $\phi'$ could be defined, which he called ``conjugate'' and for which $d\phi'=0$ in addition to $d\phi=0$.  Up to a differential form  $\omega$ this ``conjugate'' had the form of a  Grassmann dual as it was used in the general relativistic  Maxwell equation:
\beq \phi'= \epsilon_{ijkl} \frac{1}{2}\sqrt{det\, g}g^{i\alpha}g^{j\beta}B_{\alpha \beta}dx^kdx^l + c \,\omega  \, ,\footnotemark \label{eq phi conjugate}
\eeq 
where $\omega$ is independent of $\phi$ and  exact (a ``total differential''), $c$ a constant.\footnotetext{Hodge wrote this sum term by term, but used Einstein summation convention at other places of the same publication).}

Moreover,  also  $\omega$ is a harmonic form of the first kind,  and  
 as $\phi$ and $\phi'$ are the real and imaginary parts of complex differentials of the first kind, of which there are $p_g$. This and the  main theorem  showed that the second Betti number $p_2$ is given by
 \beq   p_2 = 2p_g+1 \, . \label{eq geometric genus topological}
 \eeq 
  The topological nature of the geometric genus for surfaces was thus already visible. 
 
Hodge did not stop at this point. 
Using harmonic forms he showed that the intersection matrix $A$ of  a base of the second homology of the surface is non-degenerate symmetric and the number $q$ of negative elements in the signature $\mathrm{sig}\,A = (p,q)$ (its negative index of inertia) is also $p_2$ and thus can be expressed in terms of  the geometrical genus by $q=2p_g+1$. 
The upshot of the argument was that in even stronger sense than by (\ref{eq geometric genus topological})
``\ldots $p_g$ is expressed as a topological invariant of the manifold'',
namely by the signature of its intersection matrix \citep[p. 318]{Hodge:1933genus}.

The result may appear unspectacular today; at the time it was not. It was a striking evidence that also for complex dimensions $n>1$ the theory of harmonic differential forms (and their integrals) promised  further insights into the connection between topological and birational invariants of algebraic varieties. In the words of Atiyah:
\begin{quote}
This was a totally unexpected result and, when published (8) [\citep{Hodge:1933genus}, ES], it created quite a stir in the world of algebraic geometers. In particular it convinced even the most sceptical of the importance of Hodge's theory, and it became justly famous as `Hodge's signature theorem'. Twenty years later it played a key role in Hirzebruch's work on the Riemann-Roch theorem and it remains one of the highlights of the theory of harmonic forms. \citep[p. 178]{Atiyah:1976Hodge}
\end{quote}

Of course the paper \citep{Hodge:1933genus}  stood not alone. It was a  whole series of papers in the first half of the 1930s, which taken together ``created the stir'' alluded to by Atiyah.

\subsection{\small Hodge's definition of harmonic forms, Hodge duality, and the Maxwell equation \label{subsection Hodge harmonic forms }}
Let us come back to the question how Hodge characterized  harmonic forms after 1933. In the resum\'ee paper \citep{Hodge:1935HarmonicIntegrals} he explained his differential geometric calculations  in more detail. He defined a harmonic $p$-form on a differentiable Riemannian manifold $M$ of dimension $m$,  as an  antisymmetric tensor $P$ with components 
$ P_{i_1\ldots i_p}$
 which satisfies   two conditions which in modernized notation are 
  \beq (I) \quad dP=0,\, \qquad \qquad (II) \quad  \nabla^j P_{i_1 \ldots i_{p-1} j} = 0\, . \label{eqs Hodge harmonic}
  \eeq 
  Hodge formulated them as  ``integrability conditions'' and remarked that they are the same as those stated 
already in \citep[eq. (1), (2)]{Hodge:1933genus}. He also added another form of the second condition  in terms of the tensor operations used in the  general relativistic literature for  the Grassmann complement and introduced a tensorial dualization denoted by  an upper asterix
\citep[p. 260]{Hodge:1935HarmonicIntegrals}. He soon noticed that the  prescription 
given here was too generous and refined  it in a follow up paper by the definition
\beq P^{\ast}_{j_1\ldots j_{n-p}} = \frac{1}{p!}\sqrt{g}\,  \epsilon_{i_1 \ldots i_p j_1 \ldots j_{n-p}} P^{i_1\ldots i_p} \, , \label{eq Hodge ast}
\eeq 
where $\sqrt{g}$ was the  abbreviation for $\sqrt{det\, g}$ used in the physics literature, index lifting was understood like in Ricci calculus,  the  $\epsilon$-symbol denoted the sign of the permutation, and the Einstein summation convention was assumed \citep[p. 485]{Hodge:1936Existence}.
 Concatenated with his own name the designation {\em Hodge $\ast$-operator}  for the dualization (\ref{eq Hodge ast}) would in the following years replace  the former Grassmann dual complement and become the generally used expression  for it.

Hodge  remarked  that also  $P^{\ast}$  is  an antisymmetric covariant tensor, although now in $(2m-p)$ components,  and  a ``total differential'', i.e.  closed. Condition (II) of (\ref{eqs Hodge harmonic}) can then be rewritten in terms of $P^{\ast}$. A differential form $P$ is thus {\em harmonic in Hodge's sense} \citep[p. 260f.]{Hodge:1935HarmonicIntegrals} if$\,$\footnote{Hodge wrote the following two identities in terms of the corresponding integrals.}  
\beq
 (I) \quad dP=0,\, \qquad \qquad (II) \quad    d P^{\ast}=0 \; .  \label{eqs I II Hodge harmonic} 
\eeq 
Moreover  here
\beq (P^{\ast})^{\ast} = (-1)^p P \, .
\eeq
Remember $dim \, M = 2m$; this simplifies the general relation $(P^{\ast})^{\ast} = (-1)^{p(n-p)} P$ for Riemannian manifolds of dimension $n$. 

Because of this relation (the word ``duality'' was not yet used in the 1930s by Hodge) it was clear that the harmonic forms of degree $p$ and of degree  $(n-p)$ stand in 1:1 correspondence, i.e. in later notation
\beq \Hs^p \cong \Hs^{n-p} \;, \qquad \mbox{by}  \quad P \; \mapsto \; P^{\ast}. \label{eq Hodge duality 0}
\eeq
In his influential book of 1941 Hodge introduced the language of $P^{\ast}$ as the ``{\em dual} of the form $P$'' \citep[p. 110ff.]{Hodge:1941}. 
The relation (\ref{eq Hodge duality 0})   would later become to be known as {\em Hodge duality}; after 1955 usually with an additional  dualization of vector spaces,
\[ \Hs^p \cong (\Hs^{n-p})^{\ast} \; ,  
\] 
 resulting from the bilinear pairing between $p$-forms $\alpha$ and $(n-p)$-forms $\beta$ given by 
\[ (\alpha,\beta) \longmapsto \int_M \alpha\wedge \beta \, .
\] 

The equations $(I)$ and $(II)$ of (\ref{eqs I II Hodge harmonic})  can easily be identified as a generalization of  the vacuum  Maxwell equation:  Equation (II) of (\ref{eqs Hodge harmonic}) generalizes  Einstein's form of the second Maxwell equation  (\ref{eq Maxwell II Einstein Nabla}) with the right hand side equal zero; equ (II) of. (\ref{eqs I II Hodge harmonic})  corresponds to Kottler's equ. (\ref{eq Maxwell Kottler}) and/or Cartan's (\ref{eq Cartan Maxwell}).
In fact Hodge used the notation of the   general relativistic literature. One may be tempted to consider this  as an indication that he defined harmonic forms with the analogy to the vacuum solution of the Maxwell equation in mind. The remark by  Atiyah  on Hodge's motivation quoted in the introduction (footnote \ref{fn Atiyah quote intro}) hints into the same direction.

On the other hand, Hodge neither said so in his publications of the 1930s nor  in his book  \citep{Hodge:1941}. The closest he came to allude to such a relation was, as far as I can see,  an  explanation of harmonic forms (``tensors'') in the book  as  ``the  analogues of the electrical intensity and magnetic induction'' \citep[p. 112]{Hodge:1941}. 
With this remark  Hodge referred to an analogy to classical (non-relativistic) electromagnetism, involving the electric field written as a 1-form $E=E_idx^i$ and the magnetic induction  written as a 2-form $B=B_{ij}dx^idx^j$  where $B_{ij} =\sqrt{det\, g}\, \epsilon_{ijk}\, B^k$ with $g=(g_{ij})$ some positive definite metric in dimension $n=3$ \citep[p. 111]{Hodge:1941}. The 2-form $B$ was thus written as a Grassmann-Hodge dual of the corresponding vector field with components  $B^i$.  Hodge remarked that in the case of   the Euclidean metric the classical equations of vector analysis , $curl\,E=0$ and $div\, B=0 $ were usually considered as giving rise to a scalar potential $\phi$ for the electric field and a vector potential $A=(A^i)$ for the magnetic induction. He emphasized that this is true only locally (in simply connected regions), while ``in the large'' only the vanishing of the exterior differential can be stated, $dE=0,\; dB=0$.\footnote{Hodge used a symbolic notation {\em sui generis}: $E\rightarrow 0, \; B \rightarrow 0$.}
Generously omitting  the difference of $E$ and $B$ he found it justified to continue:
\begin{quote}
We now define harmonic tensors to be the analogues of the electrical
intensity and magnetic induction {\em in the large}, and we are thus
led to the following definition: {\em A $p$-form $P$ is a harmonic form
if (1) it is regular everywhere on $M$, and (2) it satisfies everywhere the conditions} $P\rightarrow 0, \quad P^{\ast}\rightarrow 0$ [Hodge's notation at this time for $dP=0, \; dP^{\ast}=0$, ES]. \citep[p. 112. emph. in original]{Hodge:1941}
\end{quote}
 Again he did not  mention any relation to the relativistic Maxwell vacuum equation  or to the wave equation. This is surprising, although  we cannot exclude that he  he was aware of it.

A bit earlier, 
 before he presented this electromagnetic analogy,  Hodge explained the route he had taken  towards  generalizing the harmonicity condition $\Delta \phi = 0$ from functions $\phi$ to differential forms, and motivated his idea how to proceed from flat (Euclidean) space  to Riemannian geometry, again without mentioning Minkowski space or  pseudo-Riemannian manifolds. In the case of a  (positive definite) metric $g=(g_{ij})$ with Levi-Civita covariant derivative $\nabla$ and for  a function $u$ he proposed to  replace the condition for  the  flat Laplacian $\Delta$ 
\[ div\, grad\, \phi = \Delta \phi = 0
\]
by the Beltrami-Laplace operator $\Delta_g$ and 
  called  a function $\phi$ harmonic if 
\beq  \Delta_g\, \phi = \nabla_i \partial^i\, \phi = \frac{1}{\sqrt{det\, g}}\sum_{i=1}^n\partial_i (\sqrt{det\, g} g^{ij}\, \phi_j) =   0 \, , \label{eq Beltrami-harmonic}
\eeq 
with $\phi_j = \partial_j \phi$ \citep[p. 108]{Hodge:1941}. 
The last equality of (\ref{eq Beltrami-harmonic}) can be read as as
 \begin{quote}
   \ldots the condition that the $(n-1)$-form 
   \[ \frac{1}{(n-1)!} \sqrt{det\, g}\, g^{ij}\phi_j \epsilon_{i i_1\ldots i_{n-1}}  dx^{i_1}\ldots dx^{i_{n-1}} \]
   should be closed. This geometrical form of the condition suggests the generalisation of the notion of a harmonic function which we are seeking. \citep[p. 109]{Hodge:1941}
 \end{quote}
By this observation Hodge  motivated the introduction of $P^{\ast}$ like in (\ref{eq Hodge ast})  for  an alternating $p$-form $P$ as a step towards his definition of harmonic forms. He  now  talked explicitly about the {\em dual} of the form $P$ (p. 110) and showed that it has the  property
\beq P^{\ast \ast} = (-1)^{p(n-p)}P \, .\label{eq double ast}
\eeq
 Like in the first half of the 1930s he defined a {\em harmonic  $p$-form} $P$ to be one which is ``regular everywhere on $M$'' such  that both $P$ and $P^{\ast}$ are closed;  i.e. the two  equations of (\ref{eqs I II Hodge harmonic}) are satisfied \citep[p. 112]{Hodge:1941}.

 At the time when  Hodge developed his theory of harmonic forms it was well known that the vacuum Maxwell equation in flat (Minkowski) space leads to a wave equation $\Delta\, u_j=0$ for all components of the 4-potential $u= (u_j)$ of the Faraday tensor $F= du$.\footnote{With  Grassmann complement $ F^{\ast}$ and $dF^{\ast}=0$ one finds for $u=\sum_j u_jdx^j$ that $d\ast d\, u  =0 \longleftrightarrow \Delta u_j =0$ for all components $j$. }
 The vacuum Maxwell equation could thus  have supplied a striking example for Hodge's  definition of  a  harmonic differential forms.  The example had even already been adapted to  ``curved'' spaces in the literature of general relativity, to which Hodge referred indirectly by using the symbolism of the Ricci calculus with the Einstein sum convention  and other notational details.\footnote{In his bibliographic references Hodge did not include literature of theoretical physics. } 
If for Hodge the physics literature was anything more than a quarry for notations he must have  noticed that the Maxwell equation was a paradigmatic example for defining the harmonicity of differential forms in general. 
 But to my knowledge he never  mentioned  the  general relativistic vacuum Maxwell equation in his publications of the 1930s or in \citep{Hodge:1941}. This  sheds some doubts on the historical reliability of Atiyah's remark that the Maxwell equation served as a {\em motivation} for Hodge's harmonic forms. But we can also not exclude that Hodge avoided to discuss  such a relation in the written work, because for an open motivational argument two technical difficulties would have come to the fore and had, perhaps, to be discussed. First the relativistic Maxwell equation assumes a metric with Lorentzian signature; the resulting Beltrami-Laplace operator thus turns into a Beltrami-d'Alembertian and accordingly the partial differential equation becomes hyperbolic rather than elliptic. Hodge had good reasons to stay with the elliptic case. 
 
 Secondly there is a quite complicated  interrelation between the Beltrami-Laplace  characterization of the  harmonicity condition for components of differential forms and Hodge's definition. Already  a closer inspection of the Maxwell equations in a (pseudo-) Riemannian space with Levi-Civita derivative $\nabla$ and Ricci curvature $Ric = (R_{ij})$ shows for a 1-form $u=u_jdx^j$ like the 4-potential of the Maxwell field $F=du$ that the differential geometric equivalent of $dF^{\ast}= d\ast d\, u =0$ involves a curvature term. In fact it is: 
\[ \nabla_j\nabla^j\, u _i - R_{ij}u^j = 0 \qquad \longleftrightarrow  \qquad  \Delta_g  u _i -\nabla_j( \Gamma^j_{ik}u^k) - R_{ij}u^j = 0 \; 
\]  
In Hodge's approach, like in the general relativistic  Maxwell theory,  the Beltrami-Laplace harmonicity condition is thus  ``deformed'' by  additional terms depending on the Ricci curvature, the Levi-Civita connection and its derivatives.\footnote{See, e.g.,  \citep[p. 370]{Frankel:Geometry}.} 

Hodge's approach was perfectly designed to avoid the analysis of such complications. If he  was aware of this  difficulty, the  analogue to the potential of classical electromagnetism given in \citep{Hodge:1941} may perhaps be read as the expression of a didactical (over-) simplification. However,  another remark  of Atiyah on his mathematical teacher   speaks against such an  interpretation: 
\begin{quote} 
In fact Hodge knew little of the relevant analysis, no Riemannian geometry, and only a modicum of physics. His insight came entirely from algebraic geometry, where many other factors enter to complicate the picture.''  \citep[p. 186]{Atiyah:1976Hodge} 
\end{quote}
This characterization agrees  perfectly well with the discussed source texts and stands in a certain tension to  Atiyah's statement on the motivational role of the Maxwell equation for Hodge.

\subsection{\small Short remarks on the further development of Hodge's theory \label{subsection Hodge theory}}
The proof of Hodge's main theorem (\ref{eq Hodge's thm}) remained a problem for more than a decade. His first approach in \citep{Hodge:1932/34} left  gaps which he tried to fill according to a hint of H. Kneser \citep{Hodge:1936Existence}. The  second proof found wider readership when it was reproduced  in his book  \citep[chap. 3]{Hodge:1941}. It was discussed in a Princeton seminar and H.F. Bohnenblust noticed that even  the improved version was based on a problematic limit argument.\footnote{Bohnenblust gave a counter example to the limit argument,  documented in \citep[p. 1]{Weyl:1943}. }
 But shortly later,  ``building on the formal foundations laid by Hodge'', {\em Hermann Weyl} showed how the problem can be fixed  \citep{Weyl:1943}. Independently {\em Kunihiko Kodaira} developed  a proof of Hodge's theorem in his PhD dissertation, using orthogonal decomposition of the space of $p$-forms (see the passage on de Rham below). After the war  it was published in the {\em Annals of Mathematics} \citep{Kodaira:1949}  and brought him an invitation to the United States.
 
Weyl reformulated  Hodge's main theorem: For any  any differential form $f$ there exists  a uniquely determined form  $\eta \sim f$, i.e.  homologously equivalent in the sense of de Rham, for  which $ d\eta^{\ast}=0$. If  $f$ is closed, $\eta$ is  harmonic in the sense of Hodge. Weyl commented:
 \begin{quote}
 The new proposition shows at once that for any rank $p$ the space of closed forms modulo null may be identified with the space of harmonic forms. \citep[p. 6]{Weyl:1943}
 \end{quote}
In the slightly later terminology  and notation Hodge's main theorem can thus be stated as the fact, that  every cohomology class $\overline{f}$ of the de Rham cohomology $H_{dR}^k(M,\R)$ (for an at least twice differentiable manifold $M$) can be  uniquely represented by a harmonic form, or, if  $\Hs^k(M)$ denotes   the vector space of  harmonic $k$-forms on $M$,\footnote{According to \citet[p. 581]{Massey:cohomology} de Rham's theory was  explicitly formulated as a cohomology theory by H. Cartan as late as  1948 in seminars at Harvard and in his own Paris {\em S\'eminaire}. } 
\beq  \Hs^k(M) \cong H_{dR}^k(M,\R) \, . \label{eq Weyl's Hodge thm}
\eeq

Only a few years after Hodge's book appeared, {\em Kunihiko Kodaira} (1944)  and independently  {\em Georges de Rham} coauthored by {\em Pierre Bidal} (1946) introduced a Laplacian operator $\Delta_H$ adapted to the framework of Hodge's theory. Both publications introduced a new operation for $p$-forms $\omega$, denoted $\delta$ by Bidal/de Rham  and called {\em codifferential},\footnote{Kodaira used different notations and language and a divergence expression similar to Hodge's  (\ref{eqs Hodge harmonic}) in place of $d\ast$.}
\beq \delta \omega =   \ast d\hspace*{-0.2em} \ast \, \omega = (d \omega^{\ast})^{\ast} \; . \label{eq codifferential}
\eeq
In a manifold of  dimension $n$ they defined the Laplace-Hodge operator as
\beq \Delta_H\  =  (-1)^{(p+1)n} d \delta +  (-1)^{pn}  \delta d  \; .\footnotemark \label{eq Delta_H}
\eeq 
\footnotetext{\citep[p. 11]{Bidal/deRham} and up to sign \citep[p. 193]{Kodaira:1944}.}
The Swiss authors  explained that this definition can be understood 
as  a covariant generalization of the Laplacian of a vector field $A=(A^i)$ which has been used in Euclidean space in classical  electromagnetism since the late 19th century,  
\[ \Delta\, A = grad\, div \, A \; - curl\, curl\,A \; .\footnotemark
\]
 In fact, after replacing the vector field  by the corresponding  1-form $\alpha=A_idx^i$ one finds that  $\delta\alpha$ corresponds to $div\, A$ and similarly  $d\delta\, \alpha=grad\, div\, A$, $\delta d\, \alpha= curl\, curl\, A$ \citep[p. 11]{Bidal/deRham}.
 \footnotetext{The components of $\Delta\,A$ are just $\Delta\,A_i= \sum_j \partial_j^2 A_i $}
 
 Both (groups of) authors showed  that a form $\omega$ on a closed manifold $M$  (compact without boundary) satisfying
$ \Delta_H \, \omega = 0 $ 
is harmonic in the sense of Hodge.\footnote{\citep[Thm. 7]{Kodaira:1944}, \citep[p. 12]{Bidal/deRham}.}
Bidal/de Rham gave an  elegant proof  by introducing a scalar product $(\; , \;)$ on the vector space of  $p$-forms 
\[ (\alpha, \beta)=  \int_M \alpha \beta^{\ast} \, .
\]
Then $d$ and $\delta$ turned out to be adjoint operators, $(d\alpha,\beta)= (\alpha,\delta\beta)$ and vice versa. 
They found and exploited mutual orthogonality relations between forms of the types harmonic ($\Delta_H\, \alpha = 0 $), homologue zero ($\alpha= d\omega$ for some $(p-1)$-form $\omega$), ``cohomologue'' zero ($  \alpha= \delta\,\varphi $ for some $(p+1)$-form $\varphi$) and showed that every $C^2$ differential form can be decomposed in three summands, $\alpha=\alpha_1+\alpha_2+\alpha_3$, which are respectively homolog zero, cohomolog zero, and harmonic.\footnote{Because of  orthogonality of the summands $\delta d\, \alpha$ and $d\delta\,\alpha$, the harmonicity of $\alpha$ implies that both summands are zero (and vice versa). Moreover as  $(d\alpha,d\alpha)= \pm (\alpha,\delta d\alpha)$, the vanishing of $\delta d \alpha$ implies $d\alpha=0$;  similarly  $\delta \alpha=0$ \citep[p. 12]{Bidal/deRham}. } The whole article was written in an impressingly clear language in the style of the  Bourbaki group, applied to the geometry on differential manifolds and  elliptic operators. It would go beyond  the scope of the present paper, however, to discuss it in  more detail.

About the same time {\em Andr\'e Weil} opened the  study of Hodge theory on complex Hermitian manifolds with a positive definite metric given by 
\[ds^2 =\sum_{\nu=1}^n \omega_{\nu}\overline{\omega}_{\nu} = h_{\alpha {\beta}}\, dz^{\alpha} d\overline{z}^{\beta} \; \qquad \quad  (h_{\beta\alpha}=\overline{h}_{\alpha\beta})\, ,
\]
where $\nu, \alpha, \beta$ run between 1 and $n$ and the $\omega_{\nu}$ are $n$  independent linear combinations of the complex coordinate differentials $dz_1, \ldots, dz_n$ \citep{Weil:1947}.\footnote{
 A crucial paper for the study of Hermitian manifolds, cited by Weil, was \citep{Chern:1946}.} 
By demanding that the  2-form associated to the metric
\[ \omega=\sum_{\nu}\omega_{\nu}\wedge \overline{\omega}_{\nu} = \sum  \omega_{\alpha\beta}\, dz_{\alpha}\wedge\zbar_{\beta} 
\quad \mbox{with} \quad 
\omega_{\alpha \beta}=\textsf{i}\,h_{\alpha\beta} \;
 \qquad \quad (\omega_{\alpha \beta}=-\omega_{\beta\alpha})\, ,\]
 is closed, $d\omega=0$,
Weil specialized  to a K\"ahlerian metric. This allowed him to introduce a $\ast$-operator on complex differential forms  and also the codifferential $\delta$ and Laplacian $\Delta$ like in \citep{Bidal/deRham} without going back to the real structure. His main interest was directed towards extending the study of  holomorpic differential forms to meromorphic ones, i.e. those with poles. 
He proposed to generalize what Hodge had started to do with holomorphic forms on complex algebraic manifolds embedded in a projective space and endowed with a Riemannian metric (see above and \citep[p. 188ff.]{Hodge:1941}) to complex analytic manifold with K\"ahlerian metric and meromorphic forms. His proposals were soon taken up by Eckmann and Guggenheimer in a series of notes in the {\em Comptes Rendus} and continued by Kodaira.

It did not take long that also Hodge took up the thread. In \citep{Hodge:1951Kaehler} and  \citep{Hodge:1951diff-forms} he undertook a systematic study of the holomorphic and antiholomorphic forms and of mixed type $(p,q)$:
\[ \omega= \omega_{\alpha_1\ldots \alpha_p, \beta_1 \ldots \beta_q} dz^{\alpha_1}\wedge \ldots \wedge dz^{\alpha_p}\wedge d\zbar^{\beta_1}\wedge \ldots \wedge d\zbar^{\beta_q} \; ,
\] 
where $d\zbar^{\beta}$ denotes the differential with regard to the complex conjugate of the coordinate $z^{\beta}$.
The exterior differential could then be decomposed in a contribution $\partial$ of the derivative with regard to the holomorphic differentials and $\overline{\partial}$ with regard to the antiholomorphic differentials,  $d=\partial + \overline{\partial}$ and similarly for the codifferential $\delta$. This allowed to introduce harmonic forms of type $(p,q)$ and the decomposition of the harmonic forms of rank $k$ into harmonic forms of mixed type $(p,k-p)$ \citep[p. 106]{Hodge:1951Kaehler}. Let us denote the vector spaces of the latter by $\Hs^{(p,q)}(M,\omega)$ and their (real) dimensions by  $h^{(p,q)}$.  Hodge showed that the dimensions   
depend only on the complex structure, not on the K\"ahler metric. Up to isomorphism the vector spaces  thus depend only on the manifold, $\Hs^{(p,q)}(M)$  \citep[p. 109]{Hodge:1951Kaehler}.  Moreover 
\beq \Hs^k(M)= \bigoplus_{p+q=k}\Hs^{(p,q)}(M) \, , \label{eq Hodge structure}
\eeq
The dimensions $h^{(p,q)}(M)$ became to be known as {\em Hodge numbers},  the vector space including their subdivision (\ref{eq Hodge structure}) as  the {\em Hodge structure} of $M$. Because of Hodge's main theorem    $h^k =p_k$, the corresponding Betti number, while in general $h^k\neq h^{(k,0)}$.

Finally our author  formulated a necessary condition for a compact K\"ahlerian manifold  to be analytically isomorphic to an algebraic variety and called those which satisfy it  K\"ahlerian manifolds of ``restricted type'' \citep[p. 107, 110]{Hodge:1951Kaehler}, later called {\em Hodge manifold}. Three years later \citet{Kodaira:1954} was able to prove that this condition is salso sufficient. Thus Hodge's criterion turned out to be an ``intrinsic characterization'' of algebraic varieties from the standpoint of analytic manifolds.

\section{\small Finally an outlook on Hodge duality  after 1950  \label{section after 1950}}

\subsection{\small  Hodge theory and Hodge duality become sheaf cohomological \label{Hodge theory in sheaf cohomology}}
Hodge's work had an important influence on the differential geometry and topology of manifolds, cohomology theory and, of course, the geometric theory of complex functions in several variables. The latter turned into what has been called the  {\em analytic geometry} of the 20th century. In his address 
to the 1954 International Congress of Mathematics Hermann Weyl called Hodge's theory ``one of the great landmarks in the history of our science in the present century'' because it made the fruitfulness of Riemann's  ``transcendental method'' evident \citep[p. 616]{Weyl:Address}. It  would be  overconfident 
trying  to give  a resum\'ee of the  developments resulting from it. 
Only a few glimpses into one aspect of the further developments be given here, in particular the reformulation of Hodge theory in sheaf cohmological terms.\footnote{For more  information one  may  consult \citep{Weyl:Address}, \citep[pp. 254ff., 580ff.]{Dieudonne:DiffTop},  \citep{Atiyah:1976Hodge}. The ``glimpses'' are selected from these publications. \label{fn Dieudonne et al}}

Kodaira's extension of Hodge's approach to meromorphic forms (Abelian differentials of the second and third kind) bore fruit. After having constructed  forms with prescribed periods and singularities \citet{Kodaira:1949}  proved  the analogue of the famous Riemann-Roch theorem for compact 2-dimensional Riemannian manifolds and  showed that any compact K\"ahler surface  with two algebraically independent meromorphic functions can be represented as an algebraic surface.\footnote{The classical Riemann-Roch theorem deals with compact Riemann surfaces $S$ of genus $p$.  It states a relation between the  dimension $l$ of the (complex) vector space of meromorphic functions with divisor $D$,  the dimension $m$ of meromorphic differential forms on $S$ with divisor $K-D$,  the topology of the surface encoded by the genus $p$ and the order $|D|$ of the divisor: $l-m= |D|+1 - p$.  \citet[p. 664f.]{Kodaira:1949}  proved that this relation holds for compact 2-dimensional Riemannian manifolds with Hodge numbers $h^0=2A,\; h^1=2B$. 
A short explanation of the concept of  divisor is to be found in fn \ref{fn divisor}. For Riemann see  \citep[p. 182--188]{Scholz:Diss} and the  history of the  theorem until the early 20th century \citep{Gray:Riemann-Roch,Houzel:2002}. \label{fn Riemann-Roch}}
 Five years later he extended the theory and proved the result mentioned  above  that every Hodge manifold (K\"ahlerian manifold of restricte type) is bianalytically equivalent to an algebraic submanifold of a complex projective space \citep{Kodaira:1954}.  In  collaboration with D. Spencer  he was able to show that two birational invariants of an algebraic variety  $V$, introduced by Severi and called the genera $p_a$ and $P_a$  in allusion to the arithmetical genus of curves introduced by Clebsch et al., have an underpinning in the Hodge structure  of $V$.\footnote{The equality $p_a=P_a$ was known for the dimensions $n=1, 2$; for $n=3$ \citet{Severi:1909} had sketched an incomplete proof. \citet{Zariski:1952} proved a conditional equality for even $n$, if it is true for $n-1$.} 
He and Spencer introduced the {\em arithmetic genus} $a(M)$ of a compact K\"ahlerian manifold $M$ as the alternating sum of the  numbers of holomorphic forms,
\beq a(M)= \sum_{j=0}^{n} (-1)^j h^{(j,0)}(M) \, , \label{eq arithmetic genus}
\eeq
and showed that  for 
an  $n$-dimensional singularity free algebraic variety  $M=V$ the three genera are essentially the same, $ P_a=p_a=(-1)^n (a(V)-1) \; $.  
This showed that in particular  Severi's genera are  analytic and even  birational invariants.\footnote{At first Kodaira called the alternating sum of the Hodge numbers  the  ``virtual arithmetic genus''; after the proof of the identies with  Severi's genera the attribute ``virtual'' was omitted \citep[p. 642]{Kodaira/Spencer:arithmetic-genera}.}  
In their work the authors already  used the recently introduced method of sheaf cohomology.

It is impossible to sketch here the rise of sheaf theory. But it ought to be said that the assimilation of Hodge structures to sheaf theory played an non-negligible role in its early history, although it has not yet found the corresponding attention  in the historical literature.\footnote{A partial exception is the study of the origins of sheaf theory in \citep{Chorlay:2010Cousin},  for more technical  reviews in the  Bourbaki style of history one may consult \citep{Houzel:1990,Houzel:1998}. On the relation to Hodge theory see   the passages in \citep{Dieudonne:DiffTop} indicated in fn. \ref{fn Dieudonne et al}.}
Soon after  {\em Jean Leray} and {\em Henri Cartan}  introduced sheaves in the late 1940s,  {\em Pierre Dolbeault} considered  the  cohomology   with coefficients in the sheaf of
 germs of holomorphic $p$-forms on $X$, abbreviated by $\Omega^p$.\footnote{He did not specify which construction of cohomology he referred to. Since 1936 several authors had developed different approaches to cohomology theories for general spaces $X$. The most well known were probably \v{C}ech cohomology introduced by \citep{Dowker:1937} and a cohomology theory derived from $k$-cochains defined by functions on ordered  $(k+1)$ sets of points of $X$ \citep{Spanier:1948}. Not much later Hurewicz, Dugundji and Dowker showed that  Spencer's cohomology and \v{C}ech-type cohomology lead to isomorphic homology modules \citep[p. 592]{Massey:cohomology}. } 
Denoting the resulting sheaf cohomology by  $H^q(X,\Omega^p)$ he showed that in the case of a K\"ahler manifold it coincides with Hodge's harmonic forms of mixed type  \citep{Dolbeault:1953}:
 \beq  H^q(X,\Omega^p)  \cong   \Hs^{(p,q)}(X) \label{eq Dolbeault's thm}
 \eeq 
He formulated his  theorem  more generally  for any analytic variety in which case he had  to characterize the right hand side of (\ref{eq Dolbeault's thm})) by a second sheaf theoretical cohomology denoted $H^{(p,q)}(X)$, the so-called {\em Dolbeault cohomology}.  It was  derived from  co-chains $A^{(p,q)}(X)$  on a complex analytic manifold $X$ with  coefficents in the germs of mixed holomorphic-antiholomorphic  differential forms  with distributional coefficients, so-called {\em currents} of type $(p,q)$. Here the differential operator $d$ is decomposed into its holomorphic and its antiholomorphic components, $d= \partial + \overline{\partial}$, like in  the Rham cohomological interpretation of the Hodge structure. Like in (\ref{eq Dolbeault's thm}) Dolbeault then derived:
\[   H^q(X,\Omega^p)  \cong H^{(p,q)}(X) \qquad \qquad \mbox{(Dolbeault's theorem)}
\]
We cannot discuss its derivation in more detail here, but  it became important for the later work, in particular fo Serre's duality theorem. 
Moreover it  showed that $H^q(X,\Omega^p)$ may be considered as a representation of a  generalized Hodge structure  also in the general case of a complex analytic manifold   without presupposing a K\"ahlerian metric as  auxiliary device. 
   Atiyah  emphasized this conceptual shift in his report on Hodge:
\begin{quote}
For  Riemann surfaces the complex structure defines a conformal structure and hence the Riemannian metric is not far away, but in higher dimensions this relation with conformal structures breaks down and makes Hodge's success all the more surprising. Only in the 1950s, with the introduction of sheaf theory, was an alternative and more intrinsic definition given for the Hodge numbers, namely 
\[ h^{p,q} = \dim\, H^q(X,\Omega^p)
\]
where $\Omega^p$ is the sheaf of holomorphic p-forms. \citep[p. 187f.]{Atiyah:1976Hodge}
\end{quote}

This was an important basis for the generalization of  Hodge duality by  {\em Jean Pierre Serre} for any paracompact complex analytic manifold $X$.\footnote{Paracompactness of $X$ (i.e. every open covering has a locally finite refinement) is important in this  context  for constructing \v{C}ech type cohomology theories. }
\citet{Serre:1955dualite} considered the  sheaves of germs of differential forms  $A^{(p,q)}$  of type $(p,q)$ with coefficients in differentiable functions and those with distributional coefficients,  called $K_{\ast}^{(p,q)}$, and a coboundary operator $d=\partial+\overline{\partial}$ like in de Rham cohomology.\footnote{The lower star notation $K_{\ast}^{(p,q)}$ was used by Serre only from p. 17 onward; they  expressed a duality relation to the $A^{(p,q)}$.}
All this was defined 
 not only for the  manifold $X$ itself but for any complex vector bundle $V$ over $X$.
  The resulting cohomology $ H^q(X,\Omega^p(V)) $ corresponded  to Dolbeault's generalized Hodge structure and worked  not only for holomorphic differential forms on $X$  but also for meromorphic ones with singularities encoded by a divisor $D$, respectively the line bundle $V =L(D)$ associated to it.\footnote{Already Kodaira and Spencer used a  line   bundle $L(D)$  associated to a divisor for a compact K\"ahlerian manifold for dealing with meromorphic forms \citep{Kodaira/Spencer:line-bundles}.} 
  Let the  dual vector bundle be denoted by $V^{\ast}$; it is  associated to the divisor $K-D$ with $K$ the canonical divisor class.\footnote{A {\em divisor} $D$ consists of a finite collection ${D_1, \ldots, D_l}$ of analytic subvarieties of codimension 1 in $X$, endowed with integral weights $n_j\neq 0$ and written as  $D=\sum_j n_jD_j$. The {\em degree} of the divisor is $|D|=\sum_j n_j$.  The subsets $D_j$ with $n_j<0$ encode the loci of singularities (with poles of order $\leq |n_j|$), those with $n_j>0$ zeroes of order $\leq n_j$. For any (non-zero)  meromorphic function $f$ on $X$ the zeroes and poles of $f$ define a divisor  $(f)$, called {\em canonical}. The collection of all canonical divisors is  the {\em canonical divisor class}, usually denoted by $K= \{(f)\}$. \label{fn divisor}} 
In this  general setting  $A^{(p,q)}(V)$,  $K^{(p,q)}(V)$ and $K_{\ast}^{(p,q)}(V^{\ast})$ were  infinite dimensional vector spaces with a Frechet topology. This made the analysis much more demanding than in the case of a compact manifold $X$.

 Serre  stated the following generalized {\em duality theorem}  \citep[thm. 2, p. 20]{Serre:1955dualite}:
\begin{theorem}
Let $X$ be a paracompact complex analytic manifold of dimension $n$ and $V \rightarrow X$  an analytic vector bundle with dual $V^{\ast}\rightarrow X$. If in  the construction of chain complexes $A^{(p,q)}(V)$ and $K_{\ast}^{(p,q)}(V^{\ast})$, indicated above,  the antiholomorphic coboundary operator  $\overline{\partial}$   consists of vector space homomorphisms, the topological dual of the Frechet space $H^q(X,\Omega^p(V)$ constructed from the first complex $A^{(p,q)}(V)$ is isomorphic to the cohomology  $H_{\ast}^{n-q}(X,\Omega^{n-p}(V^{\ast}))$ constructed from the second complex $K_{\ast}^{(p,q)}(V^{\ast})$. In short:
\[ (H^q(X,\Omega^p(V))^{\ast} \cong H_{\ast}^{n-q}(X,\Omega^{n-p}(V^{\ast}))
\]
\end{theorem}

The theorem was a feast in dualities. It gathered up at least three (according to how one counts even five) dualities and  intertwined them into one whole: the dual relation between $p$-forms and $(n-p)$-forms  typical for the Hodge $\ast$-operator went  together with the duality between the distributional coefficients and the $K_{\ast}^{(p,q)}$ to those of the  $A^{(p,q)}$. The dualization $V^{\ast}$ of the vector bundle $V$ demanded the  substitution of the  divisor $D$ by $K -D$;  finally the Frechet space $H^q(X,\Omega^p(V)$ had to be dualized.

It  was a long path to go from the original Hodge duality (\ref{eq Hodge duality 0}) to this theorem which is rather demanding already in its formulation.  Serre did not hesitate to show that it covered special cases relevant for the study of complex analytic manifolds. Aside from a specialization for Stein manifolds (thm. 3)\footnote{A  {\em Stein manifold} ias a complex analytic manifold which is bianalytically embeddable  in a $\C^r$. An intrinsic characterization by holomorphic separability and  holomorphic convex hulls of compact subsets is possible. Non-compact Riemann surfaces are Stein manifolds.  }
 he explicated that for 
 a compact complex analytic manifold $X$ the situation becomes  close  to what Hodge had done. Because in this case the vector spaces are finite  dimensional, the theorem specializes to 
\beq  H^q(X,\Omega^p(V)) \cong (H^{n-q}(X, \Omega^{n-p}(V^{\ast}))^{\ast} \label{eq Serre's Hodge duality}
\eeq
for all $0\leq q \leq n$ \citep[thm. 4]{Serre:1955dualite}.

In the light of Dolbeault's relation (\ref{eq Dolbeault's thm})  the specialization (\ref{eq Serre's Hodge duality}) of 
 Serre's duality theorem  was clearly a generalization of Hodge's duality (\ref{eq Hodge duality 0}) to meromorphic forms. It did not need  any recourse to a metrical structure on $X$,  which could be circumvented by using dual pairings of vector spaces. Serre accomplished  for Hodge duality what  Dolbeaut had done for the core of Hodge theory, in particular Hodge's theorem.  In this sense (\ref{eq Serre's Hodge duality})  completed the transfer of Hodge theory to sheaf cohomology.
 
It became an input for a generalization of the classical theorem of Riemann-Roch to algebraic manifolds of any dimension by F. Hirzbruch.\footnote{\citet{Hirzebruch:Habil} replaced the left hand side of the Riemann-Roch theorem (see fn. \ref{fn Riemann-Roch}) by the arithmetic genus (\ref{eq arithmetic genus}) of a vector bundle $V\rightarrow X$ (in particular $V=L(D)$ for a divisor $D$) and discovered how to express the right hand side by topological invariants of $X$ and $V$ (by a polynomial in Chern classes applied to the orientation class of $H_2(X)$).}  
Soon later both theorems, Hirzebruch's and Serre's,  were even more generalized by A. Grothendieck, but both (theorems and  authors) continued to play an important role in mathematical research of the late 20th century.

\subsection{\small From Hodge duality in physics via Yang-Mills theory back to mathematics\label{subsection in physics}}
 Hodge was not really interested in  applying  his dualization in physics. Even though he  mentioned a  reference to electrodynamics in his book on {\em The Theory and Application of Harmonic Integrals} \citep{Hodge:1941} the references to physics remained  elementary and  remained at the level of the classical Maxwell equations in 3-dimensional space with a Riemannian metric. No wonder that the  main  momentum and influence of Hodge's work  was to be felt in mathematics. But it did not remain without  repercussions in theoretical physics. 
 
 During the course of years 
 the  Hodge $\ast$-operation took  the place of what had formerly been Grassmann type complements of differential forms   in  relativistic electrodynamics, although not always with explicit reference to Hodge,  e.g. \citep[p. 108ff.]{MTW}. Then the Maxwell equations acquired  a form close to the one of  de Donder, Kottler and Cartan in the 1920s (see sec. \ref{subsection GTR}):
 \beq dF= 0 , \qquad \qquad d \ast F = [4\pi]\, \mathfrak{J} \label{eq Maxwell eqs modern}
 \eeq 
 with $F$ the Faraday tensor, $\ast F$ the corresponding Maxwell tensor
  and $\mathfrak{J}$ the current density 3-form related to the current density co-vector $\mathfrak{j}$ by $\mathfrak{J} = \ast \mathfrak{j}$, and the factor $4 \pi$ is introduced, if one makes  use of Gaussian units for electromagnetism¸.\footnote{C. Kiefer made me aware of this point. For a concise discussion see \citep[p. 366ff.]{Frankel:Geometry}. }
  This view of the Maxwell equation is  particularly illuminating for a {\em premetric} approach to  electrodynamics, which tries  to avoid the reference to   a metrical structure on spacetime as far as possible. At a foundational level  the Faraday tensor $F$,  the Maxwell tensor $H$, and their fundamental equations  $dF=0$, $d H= 4\pi \mathfrak{J}$ are  introduced separately, before different   ``constitutive'' relations between the two tensors are postulated and studied like in Kottler's approach (see sec. \ref{subsection Electrodynamics}). For  relativistic Maxwell theory a metrical structure enters through the Hodge $\ast$-operator and the constitutive relation is given by 
  \beq H= \ast F \, . \label{eq constitutive relation}
  \eeq
  
In the premetric approach this  is  only one  among different alternatives, and the question may be posed whether the  ``origin'' of the metrical structure on spacetime  can be grounded on electromagnetism rather (or in addition to)   gravity
   \citep{Hehl/Obukhov:2003}.  Einstein gravity assumes the latter; here the Hodge operator plays a crucial part in expressing a basic physical law, the constitutive relation (\ref{eq constitutive relation}) of relativistic Maxwell electrodynamics.

 Hodge's $\ast$-duality also entered gravity theory, although  only in a subordinate role, e.g. for the classification of the Weyl curvature tensor \citep[p.136ff.]{Kop/Trautman}. That seems to have happened only after the rise of  Yang-Mills theories to prominence, even for authors who emphasized a Cartan geometric approach to gravity like in \citep[p. 1096ff.]{Hehl_ea:1989local_scale}.  As far as I can see,  it played  no conceptual  
  role, comparable to the one in Cartan's considerations (section  \ref{subsection Cartan}).
  
 A major field for using  Hodge type dualization arose in physics from the growing acceptance of Yang-Mills gauge theory in the standard model of elementary particle physics from 1970s onward.  This is a story of its own which still has to be told from a historical perspective.\footnote{There are many specialized articles focussing on separate themes; for sources on the origins of gauge theory see \citep{ORaif/Straumann}, for those on gauge theories of gravitation \citep{Blagojevic/Hehl}. An accessible mathematical introduction is \citep{Nielsen:Yang-Mills},  a popular account of the overall story  \citep{Crease:Creation}.}
  In their paper (1954) {\em Chen Ning Yang} and {\em Robert Mills} \nocite{Yang/Mills:1954} proposed a field theory for strong interactions between Dirac spinor fields $\psi$ of fermionic matter, extended  by (i.e. tensorized by)  dynamical degrees of freedom in a representation space of the special unitary group $SU(2)$  (the ``isospin space'') describing different but related states of elementary particles. The  theory was  modelled after the example of electromagnetism and worked with a field potential   given by  a differential form  $\As = \As_{\mu}dx^{\mu}$ on (special relativistic) spacetime
 with values in  $\mathfrak{su}(2)$, up to equivalence under so-called gauge transformations, where we abbreviate the system of coefficients as $\As_{\mu}=  (\As_{\;\; \beta \,\mu }^{\alpha}) $. 
 
From a mathematical point of view the Yang-Mills potential $\As$ can 
be understood as a  connection $\Gamma = \Gamma(\As)$ on a $SU(2)$ principal fibre bundle and properly formed associated vector bundles used to characterize fermionic matter fields. It took some time, before this interpretation became common knowledge among mathematicians and physicists.\footnote{According to Yang  he learned from differential geometers at Stony Brooks about the geometrical interpretation of gauge fields at the end of the 1960s, but started to appreciate  it only five years later, in the mid 1970s \citep[pp.  73--75]{Yang:Papers}.}
 The field strength, called the {\em Yang-Mills field} of the potential,  is then   given by a differential 2-form $\Fs= \Fs_{\mu \nu} dx^{\mu}dx^{\nu}$ with values of $\Fs_{\mu \nu}$ in $\mathfrak{su}(2)$. It arises as a covariant exterior derivative of the potential, $\Fs = d_{\Gamma}\, \As$;   geometrically it can be understood as  the curvature of the connection $\Gamma (\As)$,
 \[ \Fs_{\mu \nu } = \partial_{\mu} \As_{\nu} - \partial_{\nu} \As_{\mu} +  [ \As_{\mu}, \As_{\nu} ] \, .
 \]
The dynamics of the Yang-Mills field is assumed to be governed by a Lagrangian  analogous to the one of   Maxwell theory. This leads to   dynamical equations for the interaction field, called {\em Yang-Mills equations}, which can be written in a form similar to the Maxwell equation,
\beq d_{\Gamma}\, \Fs = 0\; , \qquad \qquad   \delta_{\Gamma}\, F =  \ast d_{\Gamma}\hspace{-0.2em} \ast\,\Fs =   4\pi \, \Js \; ,
\eeq 
with $\delta_{\Gamma}$ the co-differential analogous to (\ref{eq codifferential}). 
In contrast to (\ref{eq Maxwell eqs modern}), the 2-form $\Fs$ and the  source  1-form  $\Js$, called the Dirac current  of the present equation,  have values in $\mathfrak{su}(2)$.

 Solutions of the vacuum Yang-Mills equations, i.e. those with  $\Js=0$, are  analogue to harmonic forms. This is the case for (generalized) Yang-Mills theories with any compact group $G$ in place of $SU(2)$.  In this sense ``the general Yang-Mills theory can be considered a `{\em non-abelian Hodge theory}' '' \citep[p. 404]{Bourguignon_ea:1982}, which brings us back from physics to mathematics.

In a 4-dimensional manifold the 2-forms $\Lambda^2$  decompose into eigenspaces $\pm1$ of the Hodge operator, $\Lambda^2 = \Lambda^+ \oplus \Lambda^- $. A vacuum Yang-Mills field $\Fs$  decomposes correspondingly into $\Fs=\Fs^+ + \Fs^-$, its so-called {\em self-dual} and {\em anti-self-dual} components. Those  $\Fs$ which consist only of a (anti-) self-dual component  have been baptised {\em instantons}. 
  Shortly after the ``November revolution'' (Pickering) \nocite{Pickering:Quarks} of the standard model, which  had shown that the  paradigm of perturbatively quantized gauge field theories promised to open up a path toward  an effective model of the elementary constitution of matter, physicists started to study  the  (anti-) self-dual solutions of Yang-Mills equations in 4-dimensional {\em Euclidean} space $\E^4$. They spoke of ``pseudoparticle solutions'' and realized that these are related to topological properties of a 3-sphere bundle over $\E^4$ \citep{Belavin/Polyakov_ea}. Mathematicians soon joined; Hodge's former PhD student {\em Michael Atiyah} was among the first.\footnote{\citep{Atiyah/Hitchin_ea:1977,Atiyah:1978ICM} and diverse other papers, documented in \citep[vol. 5]{Atiyah:Works}.}
This led to the study of module spaces of instanton solutions on 4-dimensional manifolds, which became an important subfield of differential topology of 4-manifolds in the last two decades of the 20th century.\footnote{For a review see \citep[pp. 380--390]{Nash:James}.}
The resulting development, which also established a new stage for the 
 interplay between physics and geometry, may be considered  as a second chapter of the history of Hodge theory. But this is definitely a  story of its own, to be told elsewhere.

\begin{small}
\subsubsection*{Acknowledgments:} I thank  Ralf Kr\"omer for his insistence in following the trail laid out by M. Atiyah in his remarks on the physics background for Hodge's theory.  Fr\'ederic Brechenmacher and Claus Kiefer contributed with helpful remarks on a first draft of this paper.


\end{small}
\end{document}